\documentclass{article}[12pt]
\usepackage{amsfonts}
\usepackage{amsmath}
\usepackage{amssymb}
\usepackage{amscd}
\usepackage{amsthm}
\usepackage{indentfirst}
\usepackage[hmargin=3cm,vmargin=3cm]{geometry}

\newtheorem{thm}{Theorem}
\newtheorem{prop}{Proposition}[section]
\newtheorem{lma}{Lemma}

\newtheorem{cor}{Corollary}

\theoremstyle{definition}

\newtheorem{df}{Definition}[subsection]

\theoremstyle{remark}
\newtheorem*{pf}{Proof}
\newtheorem*{pfs}{Proof (sketch)}
\newtheorem{rmk}{Remark}[subsection]
\newtheorem{example}{Example}

\newcommand{\R}{{\mathbb{R}}}
\newcommand{\Z}{{\mathbb{Z}}}
\newcommand{\C}{{\mathbb{C}}}

\newcommand{\D}{{\mathbb{D}}}
\newcommand{\HH}{{\mathbb{H}}}

\newcommand{\del}{\partial}
\newcommand{\ddel}[1]{\frac{\partial}{\partial{#1}}}
\newcommand{\delbar}{\overline{\partial}}
\newcommand{\Sum}{\Sigma}
\newcommand{\G}{\mathcal{G}}
\newcommand{\Pe}{\mathcal{P}}
\newcommand{\X}{\mathfrak{X}}
\newcommand{\J}{\mathcal{J}}
\newcommand{\A}{\mathcal{A}}
\newcommand{\K}{\mathcal{K}}
\newcommand{\cH}{\mathcal{H}}
\newcommand{\ZZ}{\mathcal{Z}}

\newcommand{\cL}{\mathcal{L}}
\newcommand{\mone}{{-1}}
\newcommand{\st}{{^s_t}}

\newcommand{\intoi}{\int_0^1}
\newcommand{\til}[1]{\widetilde{#1}}
\newcommand{\wh}[1]{\widehat{#1}}
\newcommand{\arr}[1]{\overrightarrow{#1}}
\newcommand{\paph}[1]{\{ #1 \}_{t=0}^1}

\newcommand{\codim}{\text{codim}}

\newcommand{\overbar}{\overline}

\newcommand{\om}{\omega}
\newcommand{\Om}{\Omega}

\newcommand{\trace}{\text{trace}}

\newcommand{\cS}{\mathcal{S}}
\newcommand{\cD}{\mathcal{D}}
\newcommand{\cU}{\mathcal{U}}
\newcommand{\fS}{\mathfrak{S}}
\newcommand{\fk}{\mathfrak{k}}

\newcommand{\fz}{\mathfrak{z}}

\newcommand{\rJ}{\mathrm{J}}

\def\mVol{\text{Vol}(M,\omega^n)}

\begin{document}

\title{\textbf{The Action homomorphism, quasimorphisms and moment maps on the space of compatible almost complex structures}}

\author{\textsc{Egor Shelukhin}
\\
School of Mathematical Sciences \\
Tel Aviv University\\
69978 Tel Aviv, Israel \\
\tt egorshel@post.tau.ac.il}

\date{}
\maketitle

\begin{abstract}
We extend the definition of Weinstein's Action homomorphism to Hamiltonian actions with equivariant moment maps of (possibly infinite-dimensional) Lie groups on symplectic manifolds, and show that under conditions including a uniform bound on the symplectic areas of geodesic triangles the resulting homomorphism extends to a quasimorphism on the universal cover of the group. We apply these principles to finite dimensional Hermitian Lie groups like the linear symplectic group, reinterpreting the Guichardet-Wigner quasimorphisms, and to the infinite dimensional groups of Hamiltonian diffeomorphisms of closed symplectic manifolds, that act on the space of compatible almost complex structures with an equivariant moment map given by the theory of Donaldson and Fujiki. We show that the quasimorphism on the universal cover of the Hamiltonian group obtained in the second case is symplectically conjugation-invariant and compute its restrictions to the fundamental group via a homomorphism introduced by Lalonde-McDuff-Polterovich, answering a question of Polterovich; to the subgroup Hamiltonian biholomorphisms via the Futaki invariant; and to subgroups of diffeomorphisms supported in an embedded ball via the Barge-Ghys average Maslov quasimorphism, the Calabi homomorphism and the average Hermitian scalar curvature. We show that when the first Chern class vanishes this quasimorphism is proportional to a quasimorphism of Entov and when the symplectic manifold is monotone, it is proportional to a quasimorphism due to Py. As an application we show that a Sobolev distance on the universal cover of the Hamiltonian group is unbounded, similarly to the results of Eliashberg-Ratiu.

\end{abstract}

\tableofcontents

\section{Introduction and main results}

\bigskip

\subsection{Introduction}

In \cite{BargeGhysEulerBornee} Barge and Ghys have introduced a quasimorphism on the fundamental groups $\Gamma$ of surfaces of genus $g \geq 2$ (cf. \cite{PolterovichDynamicsGroups}). Their construction uses in a fundamental way the discrete action of $\Gamma$ by isometries on the hyperbolic upper half-space $\HH$. Indeed, choosing a $\Gamma$-invariant one-form $\alpha$ on $\HH$ whose differential is bounded in the way $|d\alpha|\leq C_\alpha |\sigma_\HH|$ for a constant $C_\alpha$ with respect to the hyperbolic Kahler form $\sigma_\HH$ on $\HH$, the quasimorhism is given by integrating $\alpha$ over the geodesic $l(x,\gamma\cdot x)$ between a fixed base-point $x$ and its image $\gamma\cdot x$ under the action of an element $\gamma \in \Gamma$. Using these quasimorphisms Barge and Ghys have obtained results on the second bounded cohomology $H_b^2(\Gamma)$ of such groups $\Gamma$. Further results on the second bounded cohomology of discrete groups following from their actions upon certain spaces with "negative enough" curvature - e.g. Gromov-hyperbolic groups - were studied extensively in \cite{EpsteinFujiwara,Fujiwara,HamenstadtDiscrete,Hamenstadt,MonodShalom}  to name a few works in such a direction. The second bounded cohomology of finite dimensional Lie groups was also studied extensively. For example, in the works \cite{GuichardetWigner,DupontGuichardet} and others, the action of simple Hermitian symmetric Lie groups $G$ upon their symmetric space $X = G/K$ of non-compact type was utilized to construct bounded $2$-cocycles on $G$. The basic construction of such cocycles similarly uses the integration of the natural Kahler form $\sigma_X$ on $X$ on simplices with geodesic boundaries.

We shall first formulate a general setting in terms of the action of a group $G$ on a space $X$ for constructions related to integration on geodesic simplices to yield bounded $2$-cocycles. Then we formulate a general principle, again in terms of such actions, for the construction of primitives to such cocycles in the (unbounded) group cohomology, to wit - quasimorphisms - functions that satisfy the homomorphism property up to a uniformly bounded error. For one, our construction gives a symplectic formula for the quasimorphisms on the universal covers $\til{G}$ of simple Hermitian symmetric Lie groups whose differentials equal the Guichardet-Wigner cocycles (cf. \cite{GuichardetWigner,DupontGuichardet,ClercKoufany,ShternAutomaticContinuity,BurgerIozziWienhardMaximalToledoInvariant}). A key notion in our construction is the use of \emph{equivariant moment maps} for the Hamiltonian action of a group $\G$ on a space $\X$ with a symplectic form $\Om$. Another key notion is that of the Action homomorphism of A.Weinstein \cite{Weinstein} that generalizes to general Hamiltonian actions with equivariant moment maps. As our construction is rather formal, or "soft" in the terminology of Gromov \cite{GromovSandH} in that it does not require the solution of partial differential equations or the convergence of certain series, it readily applies to the infinite dimensional case.

Indeed there have been many constructions of equivariant moment maps for actions of infinite dimensional Lie groups on infinite dimensional symplectic spaces $(\X,\Om)$. Starting with the work of Atiyah and Bott \cite{AtiyahBottConnections,AtiyahBottConnectionsPrelim} - for the action of gauge groups of principal bundles over Riemann surfaces on the corresponding spaces of connections, with numerous later developments including an extension to higher dimensions - a general framework for the Hitchin-Kobayashi correspondence \cite{DonaldsonComplexSurfacesHitchinKobayashi,UhlenbeckYauGeneralHitchinKobayashi,DonaldsonInfiniteDeterminantsHitchinKobayashi}, the works of Donaldson \cite{DonaldsonMomentMapsDiffGeom,DonaldsonMomentMapsDiffeomorphisms,DonaldsonGaugeTheoryComplexGeometry} and  Fujiki \cite{Fujiki} for actions of diffeomorphism groups upon spaces of mappings (submanifolds or sections of bundles), and more recent advances e.g. \cite{FutakiMomentMaps,Fine} this has been an active and fruitful area of research for over three decades, with many applications - for example to Kahler geometry. Of these the Donaldson-Fujiki \cite{DonaldsonGaugeTheoryComplexGeometry,Fujiki} framework of the scalar curvature as a moment map for the action of the Hamiltonian group on the space of compatible almost complex structures fits the setting of our construction. We shall, therefore, apply this framework to build new quasimorphisms on the Hamiltonian group, or its universal cover, of an arbitrary symplectic manifold of finite volume (and of an arbitrary closed symplectic manifold in particular). Similarly to the finite-dimensional case, our quasimorphism provides a group-cohomological primitive for the restriction to the Hamiltonian group of a certain $2$-cocycle that was constructed using the natural notion of geodesic simplices in spaces of almost complex structures by Reznikov \cite{ReznikovCocyclesVol,ReznikovCharClassSymp,ReznikovAnalyticTop} in his studies of the cohomology of the group of symplectomorphisms.

The intriguing topic of the study of quasimorphisms on groups of (Hamiltonian) symplectomorphisms has a long history. A very early work of Eugenio Calabi \cite{CalabiHomomorphism} constructs a homomorphism on the group of compactly supported symplectomorphisms of the symplectic ball of arbitrary dimension $2n$. An early example of a quasimorphism on a symplectomorphism group that is not a homomorphism was constructed by Ruelle \cite{Ruelle} on the group of compactly supported volume preserving diffeomorphisms of the two-dimensional disk, as a certain average asymptotic rotation number. This result was generalized using the Maslov quasimorphism on the universal cover $\til{Sp}(2n,\R)$ of the linear symplectic group by Barge and Ghys \cite{BargeGhysEulerMaslov} to the group of compactly supported symplectomorphisms of the symplectic ball of arbitrary dimension $2n$. A quasimorphism on the universal cover $\til{Symp}(M,\om)$ of closed symplectic manifolds $(M,\om)$ with $c_1(TM,\om) = 0$ was rather recently constructed by Entov \cite{EntovCommutatorLength}, generalizing the previous quasimorphism in the sense that it equals the Barge-Ghys average Maslov quasimorphism when restricted to each subgroup of diffeomorphisms supported in an embedded ball - we shall say that it has the Maslov local type. In a recent work of Py \cite{pyqm,PyThesis} a quasimorphism on $\til{Ham}(M,\om)$ for closed symplectic manifolds $(M,\om)$ with $c_1(TM,\om) = \kappa [\om]$ for $\kappa \neq 0$ was constructed as a rotation number using the notion of a prequantization of an integral symplectic manifold. The local type of the Py quasimorphism is Calabi-Maslov - it equals a certain linear combination of the Calabi homomorphism and the Barge-Ghys average Maslov quasimorphism when evaluated on diffeomorphisms supported in a given embedded ball. A compelling discovery of quasimorphisms of Calabi local type was made by Entov and Polterovich in \cite{EntovPolterovichCalabiQmQh} - one distinctive feature of which is that the embedded balls should be small enough - using "hard" methods of Hamiltonian Floer homology and the algebraic properties of quantum homology. These methods were since generalized and extended to a large class of manifolds \cite{YaronCalabiQmNonnMonotone,EntovPoletrovichSymplecticQuasistatesSemiSimplicity,YaronTyomkin,UsherDualityFloerNovikovFiltration,UsherSpectralNumbersFloer}, a very recent result due to Usher \cite{UsherDeformedQHandCalabiQm} showing e.g. the existence of Calabi quasimorphisms on $\til{Ham}$ of every one-point blowup of a closed symplectic manifold. 
The sequent question of constructing a "soft" quasimorphism of Calabi local type on the Hamiltonian group of a closed symplectic manifold was recently solved for the two-torus and for surfaces of genus $g \geq 2$ by Py \cite{pyqm}. The first case builds upon the works of Ghys and Gambaudo \cite{GhysGambaudoEnlacements,GhysGambaudoCommutatorsSurfaces} in dimension $2$ that describe the Calabi homomorphism and a large number of quasimorphisms, using such methods as the action of diffeomorphism groups upon the configuration spaces of distinct points in a surface (these works have been since developed in many other papers - cf. \cite{BrandThesis}). The second case uses prequantizations and the notion of the bounded Euler class (which is again related to the boundedness of the symplectic area of geodesic triangles), and can be extended to compact quotients of simple Hermitian symmetric spaces $X$ of non-compact type by discrete groups of isometries \cite{PyThesis}. 
Another quasimorphism on $\til{Ham}(M,\om)$ for $(M,\om)$ the complex projective space $(\C P^n, \om_{FS})$ with the natural Fubini-Study Kahler form can be derived from the work of Givental \cite{GiventalQuasimorphism} that uses methods of generating functions, which also has the Calabi property by the work of Ben Simon \cite{GabiCalabiQm} and can easily be shown to descend to $Ham(M,\om)$ itself by results from \cite{LoopRemarks}.
In fact necessary and sufficient conditions for the above quasimorphisms on a group $\til{\G}$ to descend to $\G$ are given by the vanishing of certain homomorphisms $\pi_1(\G) \to \R$. This happens automatically for sufaces where the fundamental group of $\G = Ham(M,\om)$ is finite, which is also known to be the case for certain four-dimensional symplectic manifolds - e.g. $(\C P^2,\om_{FS}),(\C P^1 \times \C P^1, \om_{FS} \oplus \om_{FS})$ \cite{GromovPseudohol} (cf. \cite{JHolSymp}). Remarkably, for all monotone examples - $(M,\om)$ such that $c_1(TM,\om) = \kappa [\om]$ for $\kappa \neq 0$ - the homomorphism is the same one \cite{RigidSubsets} - the Action-Maslov homomorphism of Polterovich \cite{PolterovichLoops} (cf. \cite{LoopRemarks}).

The quasimorphism we construct has Calabi-Maslov local type - it restricts to the difference of suitable multiples of the Calabi homomorphism \cite{CalabiHomomorphism,IntroSymp} and of the Barge Ghys average Maslov quasimorphism on the subgroup of Hamiltonian diffeomorphisms supported in a small ball. Its restriction to the fundamental group of $\G$ is equal by construction to the generalized Action homomorphism, involving in this case the Hermitian scalar curvature, and is also computed via a homomorphism earlier introduced in \cite{LMP} using a Hamiltonian fiber bundle obtained by the clutching construction. A previous work that applies the theory of the Hermitian scalar curvature as a moment map to the study of the topology of the Hamiltonian group is \cite{AbreuGranjaKitchlooSketch, AbreuGranjaKitchlooBIG}.

Furthermore, our quasimorphism agrees with the quasimorphisms of Py and Entov whenever these quasimorphisms are defined. While, having a Maslov component in the local type, our quasimorphism can at best be continuous in the $C^1$-topology, it is rather easily seen to be coarse-Lipschitz in the Sobolev $L^2_2$-metric, using the isoperimetric propery of Kahler manifolds with a bounded primitive of the Kahler form. This allows us to prove that the Sobolev $L^2_2$-metric is unbounded on $\til{\G}$ of every symplectic manifold of finite volume, extending a consequence from previous works of Eliashberg-Ratiu \cite{EliashbergRatiu} on the $L^2_1$-metric in the case when the symplectic manifold is exact. Moreover, we show that on manifolds like the blowup $Bl_1(\C P^2)$, where the restriction of the quasimorphism to $\pi_1 \G$ does not vanish, the metric is not bounded on $\pi_1 \G$ either. 
We conclude with some questions and discussion related to the topics presented in the paper.

As an aside, it is curious to note that this paper touches upon two directions that both have their origins with Eugenio Calabi - the study of canonical metrics on Kahler manifolds (e.g. \cite{CalabiCY,CalabiExtremalMetrics1,CalabiExtremalMetrics2}) and the theory of the Calabi homomorphism (\cite{CalabiHomomorphism}).

\subsection{Moment maps}\label{section - Moment maps}
Assume that a Lie group $\mathcal{G}$ acts $\G \times \mathfrak{X} \to \mathfrak{X}, (g,x) \mapsto g \cdot x$ on a symplectic manifold $(\mathfrak{X},\Omega)$ in a Hamiltonian fashion. Here both the group and the manifold can be infinite dimensional. The action gives a homomorphism $\G \to Diff(\X)$, $\phi \mapsto \overline{\phi}$ with the property that to each element $X \in Lie (\G)$ there corresponds an element $\mu(X) \in C^\infty(\X,\R)$, such that

\begin{enumerate}
\item the equation $\iota _\Xi \Om = - d \mu(X)$ holds for $\Xi \in V.F.(\X)$ - the vector field on $\X$ corresponding to $X$
\item the resulting map $Lie (\G) \to C^\infty(\X,\R)$ is a homomorphism of Lie algebras (the Lie structure on the latter is given by the Poisson bracket of the symplectic form $\Omega$).
\end{enumerate}

The second condition is equivalent to the linearity and {\it equivariance} of the map $X \mapsto \mu(X)$ - for all $X \in Lie (\G)$ and $\phi \in \G$ we have $$\mu(Ad_{\phi}X) = \mu(X)\circ\overline{\phi}\,^{-1}.$$ In one direction one differentiates this equality and the other can be found in \cite{IntroSymp} Lemma 5.16.

Note that the map $X \mapsto \mu(X)$ gives us a pairing $\mu: Lie (\G) \times \X \to \R$ that is linear in the first variable, and therefore a map $x\mapsto \mu(-)(x): \X \to (Lie (\G))^*$. The equivariance condition corresponds to the invariance of the pairing with respect to the diagonal action of $\G$ - for all $X \in Lie (\G), x \in \X$ and $\phi \in \G$ we have \[\mu(Ad_\phi X)(\phi \cdot x) = \mu(X)(x).\]
We call $\mu$ in any one of these three equivalent formulations a {\it moment map} for the Hamiltonian action of $\G$ on $\X$.

\begin{rmk}
For infinite-dimensional Lie groups we use the approach of regular Fr$\acute{\text{e}}$chet Lie groups (cf. \cite{MilnorInfiniteDimensional} and references therein), while one could also use the inverse limit (ILH or ILB) approach of Omori \cite{Omori}. In any case, as we are interested only in the soft features of the theory of Lie groups and our infinite-dimensional example is a diffeomorphism group where all computations can be carried out as explicit differential-geometric formulae, the foundational theory of infinite-dimensional Lie groups can for the most part be ignored. The same remark applies to infinite-dimensional symplectic manifolds.
\end{rmk}

\subsection{The action homomorphism}\label{action}

Assume that $\pi_1(\X)=0$. Denote by $\Pe_\Om \subset \R$ the spherical period group $\langle\Om,\pi_2(\X)\rangle$ of $\Om$. Following Weinstein \cite{Weinstein}, we define the Action homomorphism $\pi_1(\G) \to \R/\Pe_\Om$ as follows.

Suppose a class $a \in \pi_1(\G)$ is represented by a path $\{\phi_t\}$ based at the identity element $Id$. Pick a point $x \in \X$. Consider its trace $\phi_x = \paph{\phi_t \cdot x}$ under the action of the loop. Pick a disk $D$ that spans $\phi_x$ - that is $D: \D \to \X$ is a smooth map from $\D=\{|z|\leq 1\} \subset \C$ to $\X$ that satisfies $D(e^{2\pi i t})=\phi_t \cdot x$ for all $t\in S^1=\R/\Z$. Then the Action homomorphism is defined as \[\mathcal{A}_\mu(a)= \int_D \Omega - \int_0^1\mu(X_t)(\phi_t \cdot x)\,dt.\]

It is independent of $x \in \X$ by the first property of $\mu$ and of $\{\phi_t\}$ in the homotopy class $a \in \pi_1(\G,Id)$ by the second property of $\mu$. It does depend on the spanning disk $D$, however the ambiguity lies in $\Pe_\Om$. At last, the homomorphism property follows by a short concatenation argument. Detailed proofs can be found in Section \ref{Proofs}.

\begin{rmk}
Note that when $\pi_2(X)=0$, the Action homomorphism takes values in $\R$, since $\Pe_\Om = 0$.
\end{rmk}

\begin{rmk}\label{new action extends original}
This definition extends the original definition because given a closed symplectic manifold $(M,\omega),$ the group $\G = Ham (M, \omega)$ acts on $(M,\omega)$ in a Hamiltonian fashion with the equivariant moment map $\mu(X) = H_X$ where $H_X \in C^\infty(M,\R)$ is the zero-mean normalized Hamiltonian function of $X$. On an open symplectic manifold $(M,\omega)$ the group $\G = Ham_c (M, \omega)$ of compactly supported Hamiltonian diffeomorphisms acts in a Hamiltonian fashion with the equivariant moment map $\mu(X) = H_X$ where $H_X$ is the compact-support normalized Hamiltonian function of $X$. To ensure the existence of a contracting disk, we assume that the manifold is simply connected in the open case. In the closed case the contracting disk always exists by Floer theory, by the existence of the Seidel element or by a direct geometric degeneration argument \cite{HamTracesContract}.
\end{rmk}

\subsection{Preliminaries on quasimorphisms}\label{preliminaries on quasimorphisms}

A \emph{quasimorphism} $\nu$ on a group $\G$ is a function $\nu: \G \to \R$ that satisfies the additivity property up to a uniformly bounded error. That is for all $x \in \G$ and $y \in \G$ we have \[\nu(xy) = \nu(x) + \nu(y) + b(x,y),\] where \[|b(x,y)| \leq C_{\nu}\] for a constant $C_{\nu}$ depending only on $\nu$ (and not on $x,y$). In such cases the limit \[\overbar{\nu}(x):= \lim_{k \to \infty} \frac{1}{k}\nu(x^k)\] exists by Fekete's lemma on subadditive sequences and is also a quasimorphism. Moreover, it is \emph{homogenous} that is \[\overbar{\nu}(x^k) = k \;\overbar{\nu}(x)\] for all $x \in \G$ and $k \in \Z$ and satisfies \[\overbar{\nu} \simeq \nu,\] where for any two functions $a, b: \G^m \to \R$ we write \begin{equation}\label{notation - equal up to bounded difference}a \simeq b\end{equation} if they differ by a uniformly bounded function $d: \G^m \to \R$ - that is $|d(x_1,...,x_m)| \leq C_d$ for a constant $C_d$ independent of $x_1,...,x_m$. We refer to the book \cite{CalegariScl} by Calegari for these statements and for additional information about quasimorphisms.

We will use the following simple fact.

\begin{lma}\label{nu of x and nu of x^{-1}}
For every quasimorphism $\nu:\G \to \R$ we have $\nu(x) \simeq - \nu(x^{-1})$ as functions $\G \to \R$.
\end{lma}

\begin{pf}
Indeed $\nu(x) \simeq \overbar{\nu}(x) = - \overbar{\nu}(x^{-1}) \simeq - \nu(x^{-1})$.
\end{pf}

Explicit constructions of quasimorphisms on Lie groups often use rotation numbers. For this purpose we require the notion of the variation of angle of a continuous path $\delta:[0,1] \to S^1.$

\begin{df}\label{varangle definition}
We define the full \emph{variation of angle} of $\delta:[0,1] \to S^1$ as \[varangle(\delta)=\til{\delta}(1) - \til{\delta}(0)\] for any continuous lift $\til{\delta}:[0,1] \to \R$ of $\delta$ to the universal cover $\R \xrightarrow{\Z} S^1$.
\end{df}

\subsection{A general principle for constructing quasimorphisms}\label{quasimorphisms}
The general principle says that when groups act well enough on spaces of negative enough curvature, then
they have quasimorphisms and non-trivial bounded (or bounded-continuous) cohomology. While usually this principle is applied to proper discontinuous actions of discrete groups, we propose a version of this principle for smooth actions of (possibly infinite dimensional) Lie groups. Firstly, we  propose  a version of "negative enough curvature" - (possibly infinite-dimensional) symplectic manifolds $(\X,\Om)$ with bounded Gromov norm of $\Om$. We make, more specifically, the following definition.

\begin{df}(Domic-Toledo space $(\X,\Om,\K)$)
Assume that $\X$ has $\pi_1(\X) = 0$ (as before) and $\pi_2(\X) = 0$ also. Moreover assume that there is a system $\mathcal{K}$ of paths $[x,y]:=\gamma(x,y)$ for all $x \in \X$ and $y \in \X$, such that for all $x,y,z \in \X$ $$|\int_{\Delta(x,y,z)} \Om| < C_\X,$$ for a constant $C_\X$ that does not depend on $x,y,z$. Here $\Delta = \Delta(x,y,z)$, which we will call a \emph{geodesic triangle} is any disk with boundary $\partial \Delta = [x,y]\cup [y,z] \cup [z,x]$. We call the triple $(\X,\Om,\K)$ a \emph{Domic-Toledo} space.
\end{df}

Next we propose a version for "act well enough" - by "isometries" with an equivariant moment map. More exactly, we make the following definition.

\begin{df}(Hamiltonian-Hermitian group $\G$) We call a (possibly infinite dimensional) Lie group $\G$ \emph{Hamiltonian-Hermitian} if it acts on a Domic-Toledo space $(\X,\Om,\K)$ - preserving $\K$ and $\Om$ - with an \emph{equivariant moment map} \[\mu:\X \times Lie(\G) \to \R.\] We say that the action of $\G$ on $(\X,\Om,\K)$ preserves $\K$ if for every two points $x \in \X$ and $y \in \X$ and every $g \in \G$ we have \[g\cdot[x,y] = [g\cdot x,g \cdot y].\]
\end{df}

\begin{rmk}
All examples of Domic-Toledo spaces known to the author are (possibly infinite-dimensional) Kahler manifolds $(\X,\Om,J)$ with $[x,y]$ being the geodesic segment between $x \in \X$ and $y \in \X$. A first set of examples is given by Hermitian symmetric spaces $\mathcal{D}$ of non-compact type (bounded Hermitian domains) \cite{DomicToledo,ClercOrsted}. The second one (trivially containing the first) is given by spaces of global sections of bundles with fiber $\cD$ over a manifold $(M,\phi)$ with a volume form $\phi$ of finite volume.
\end{rmk}

\begin{rmk}
Examples of finite dimensional Hamiltonian-Hermitian groups are given by Hermitian symmetric Lie groups - like $Sp(2n,\R)$  -  since they act by Hamiltonian biholomorphisms on the corresponding symmetric spaces of non compact type equipped with the Bergman Kahler structure, which is Kahler-Einstein. Therefore, the natural lift (by use of the differential) of these diffeomorphisms to the top exterior power of the tangent bundle furnishes the action with an equivariant moment map (note that the Kahler-Einstein condition implies that ($-i$) times the curvature of the Chern connection on these bundles, given by the Hermitian metric, equals to the Kahler form on one hand, and on the other hand the corresponding connection form is surely preserved by the lifts). Details are presented in Section \ref{FinDimExamples}.

Infinite dimensional examples are given by groups $Ham(M,\om)$ of closed symplectic manifolds $(M,\om)$ since these act of the spaces $\J$ of compatible almost complex structures, which is a Domic-Toledo space - since it is the space of global sections of a bundle over $(M,\om)$ with fiber the Siegel upper half-plane. This class of examples can be extended to arbitrary symplectic manifolds of finite volume. Details are presented in Section \ref{quasimorphisms on the Hamiltonian groups of symplectic manifolds}.
\end{rmk}

We now construct a quasimorphism on the universal cover of a Hamiltonian-Hermitian group $\G$ with an equivariant moment map $\mu$ and Domic-Toledo space $(\X,\Om, \K)$. Given a path $\paph{g_t}$ in $\G$ with $g_0 = Id$, $g_1 = g$ representing a class $\til{g}$ in $\til{\G}$, consider the loop $\paph{g_t\cdot x} \# \;[g\cdot x, x]$ for a fixed basepoint $x \in \X$. Fill it by any disk $D = D_{\paph{g_t}}$. Then define \begin{equation}\label{nu_x equation}\nu_{x}(\til{g})= \int_{D}\Om - \int_0^1\mu(X_t)(g_t \cdot x)dt,\end{equation} where $\paph{X_t}$ is the path in $Lie(\G)$ corresponding to the path $\paph{g_t}$. In Section \ref{proof of general principle for Hamiltonian-Hermitian} we show that this value is well-defined and gives a real-valued quasimorphism $\nu_x:\til{\G} \to \R$ on the universal cover of $\G$.

\begin{thm}\label{General principle for Hamiltonian-Hermitian}
Any Hamiltonian-Hermitian group $\G$ acting with an equivariant moment map $\mu$ on the corresponding Domic-Toledo space $(\X,\Om, \K)$ admits a real-valued quasimorphism $\nu_x:\til{\G} \to \R$ on its universal cover for each point $x \in \X$, given by Equation \ref{nu_x equation}. Moreover, the homogeneization $\nu$ of $\nu_{x}$ does not depend on the basepoint $x$. By construction, the quasimorphism $\nu$ restricts to the homomorphism $\A_\mu$ on $\pi_1(\G)$.
\end{thm}

\begin{rmk}
If we assume additionally that the loop $[x,x] \in K$ is the constant path at $x$, then the quasimorphism $\nu_x$ also restricts to $\A_\mu$ on $\pi_1(\G)$.
\end{rmk}

Note that this theorem does not state that the homogenous quasimorphism $\nu$ is necessarily not a homomorphism, or even not trivial. It can in principle be identically equal to zero. However, in all the known examples it turns out to be non-trivial and not a homomorphism.

The key feature of the proof which we defer to Section \ref{proof of general principle for Hamiltonian-Hermitian} is that the differential of $\nu_x$ in group cohomology satisfies \begin{equation}\label{differential of nu_x as integral over geodesic simplex}b(g,h)=\nu_{x}(\til{g}\til{h}) -\nu_{x}(\til{g}) - \nu_{x}(\til{h}) = \int_{\Delta(x, g\cdot x, gh \cdot x)} \Om\end{equation} for $\til{g},\til{h} \in \til{\G}$ with endpoints $g,h \in \G$. The latter is a bounded cocycle by the properties of Domic-Toledo spaces and "isometric" actions upon them.

\begin{rmk}\label{nu_x and -nu_x}
From Equation \ref{differential of nu_x as integral over geodesic simplex}, given that for all $x \in \X$, $[x,x]$ is the constant path at $x$, it follows that for all $\til{\phi} \in \til{G}$ we have \[\nu_x(\til{\phi}^{-1}) = -\nu_x(\til{\phi}).\] Indeed the difference equals $\int_{\Delta(x,\phi\cdot x, x)} \Om = 0$, since we can choose a degenerate filling disk.
\end{rmk}

Furthermore, we would like to explore the invariance of the quasimorphism with respect to larger groups extending a given action of a Hamiltonian-Hermitian group $\G$ on a Domic-Toledo space. For this we have the following proposition, which we prove in Section \ref{proof of general principle for Hamiltonian-Hermitian}.

\begin{prop}\label{general symp independence}
Assume $\G \subset \mathcal{H}$ is a normal subgroup, $\G$ is a Hamiltonian-Hermitian group acting with an equiavariant moment map $\mu$ on the Domic-Toledo space $(\X,\Om,\K)$, and $\mathcal{H}$ is a (possibly infinite-dimensional) Lie group that acts on $(\X,\Om,\K)$ preserving $\Om$ and $\K$ and extending the action of $\G$ (however not necessarily with a moment map).  Assume moreover, that the moment map $\mu: Lie(\G)\times \X \to \R$ is equivariant with respect to the action of $\cH$ (note that as $\G \subset \cH$ is normal, $\cH$ acts on $Lie(\G)$ by the adjoint representation). Then $\nu_x(h\til{g}h^{-1})=\nu_{h^{-1}x}(\til{g})$ for all $\til{g} \in \til{\G}$ and $h \in \mathcal{H}$. Consequently, by the independence of the homogeneization upon the basepoint, we have \[\nu(h\til{g}h^{-1})=\nu(\til{g}),\] for all $\til{g} \in \til{\G}$ and $h \in \mathcal{H}$. Equivalently $\nu(\til{h}\til{g}\til{h}^{-1})=\nu(\til{g}),$ for all $\til{g} \in \til{\G}$ and $\til{h} \in \til{\mathcal{H}}$.
\end{prop}

\subsection{The scalar curvature as a moment map}\label{scalar}

Given a compact symplectic manifold $(M,\omega)$ consider the space $\mathcal{J}$ of $\omega$-compatible almost complex structures.
This space can be given the structure of an infinite-dimensional Kahler manifold $(\J,\Om,\mathbb{J})$ as follows. Consider the bundle $\cS \to M$, the general fibre of which over $x \in M$ is the space $\rJ_c(T_x M, \om_x) \cong Sp(2n)/U(n)$ of $\om_x$-compatible complex structures on $T_x M$. As $\rJ_c$  posses a canonical $Sp(2n)$-invariant Kahler form $\sigma = \sigma_{trace}$ we have a fiberwise-Kahler form $\sigma$ on $\cS$. Note now that $\J = \Gamma(M,\cS)$ - the space of global sections of the bundle $\cS \to M$. Now define $\Om(A,B) := \int_M \sigma_x(A_x,B_x) \om^n(x)$. The complex structure $\mathbb{J}$ on $\J$ is defined as $\mathbb{J}_J A = JA$ for $A \in T_J \J$. Surely $\Om$ and $\mathbb{J}$ are compatible.

Note that the group $\G = Ham (M, \omega)$ of Hamiltonian diffeomorphisms acts on $\J$ by $\phi \cdot J := \phi_*J$.
This action can be shown to be Hamiltonian \cite{DonaldsonGaugeTheoryComplexGeometry,Fujiki} with respect to the form $\Om$.  The moment map is given as follows.

First note that the Lie algebra of $\G$ is isomorphic to the space $C^\infty(M,\R)/\R \cong C^\infty_0(M,\R)$. The latter space consists of smooth functions $F$ on $M$ with integral zero: $\int_M F \om^n = 0$. For an element $\phi \in \G$, the adjoint action is given in these conventions by \begin{equation}\label{Adjoint action in G on C^infty}Ad_\phi H = (\phi^{-1})^*H.\end{equation}

To a function $H \in Lie (\G) \cong C^\infty_0(M,\R)$ there corresponds the function $\mu(H)$ on $\J$ given \cite{DonaldsonGaugeTheoryComplexGeometry,Fujiki} by the formula: \begin{equation}\label{definition of mu for G=Ham}\mu(H)(J) = \int_M S(J)H \om^n,\end{equation} where $S(J) \in C^\infty(M,\R)$ is the \textit{Hermitian scalar curvature} of the Hermitian metric $h(J)= g(J) - i \omega $ defined as follows. Consider the Hermitian line bundle $L=\Lambda_\C^n(TM,J,h(J))$. It has a natural connection $\nabla^n$ induced from the canonical connection $\nabla$ on $(TM,J,h(J))$ (cf. \cite{Gauduchon} Section 2.6, \cite{KobayashiConnectionsAlmsotComplex} and \cite{TosattiWeinkoveYauTamingSymplecticForms} Section 2, and references therein) defined by the properties \[\nabla J = 0, \nabla h = 0, T^{(1,1)}_\nabla = 0.\] This connection can also be equivalently (by \cite{Gauduchon} Section 2) defined by use of $\delbar$-operators, as in \cite{DonaldsonGaugeTheoryComplexGeometry}. The connection $\nabla^n$ has curvature $i\til{\rho}$ for the lift $\til{\rho}$ of a real valued closed two form $\rho \in \Om^2(M,\R)$ on $M$ by the natural projection $L \to M$. We define $S(J) \in C^{\infty}(M,\R)$ by \begin{equation}\label{definition of S(J)} S(J)\om^n = n \rho \wedge \om^{n-1}.\end{equation} Whenever $J$ is integrable $S(J)$ coincides with the scalar curvature of the Riemannian metric $g(J)$. In the above $g(J)$ is the Riemannian metric corresponding to $J$ given by $g(J)(\xi,\eta)=\om(\xi,J\eta)$. Note that $g(\phi_*J) = (\phi^{-1})^*g(J)$ and consequently the same is true for $h(J)$. Hence \begin{equation}\label{}S(\phi_*J) = (\phi^{-1})^*S(J)\end{equation} for all $\phi \in \G$.

From the equalities \ref{Adjoint action in G on C^infty},\ref{definition of mu for G=Ham} and \ref{definition of S(J)} we obtain $\mu(H)(\phi^{-1}_*J) = \int_M S(\phi^{-1}_*J)H \om^n = \int_M \phi^*S(J)H \om^n = \int_M S(J) (\phi^{-1})^*H \om^n = \mu((\phi^{-1})^*H)(J) = \mu(Ad_\phi H)(J)$. Therefore the moment map is equivariant.

We remark that the action of $\G$ on $\J$ can be extended to the action of $\cH = Symp(M,\om)$ that preserves $\Om$ and $\J$. Moreover $\G \subset \cH$ is a normal subgroup and (by the same computation as above) the moment for map $\mu: Lie(\G) \times \J \to \R$ for the action of $\G$ on $\J$ is equivariant with respect to the action of $\cH$ (which acts on $Lie(\G)$ by the adjoint action $Ad_\psi H = (\psi^{-1})^*H$, $\psi \in \cH$).


\subsection{Quasimorphisms on the Hamiltonian groups of symplectic manifolds}\label{quasimorphisms on the Hamiltonian groups of symplectic manifolds}

Here we apply the general principle for constructing quasimorphisms to the group $\G = Ham(M,\om)$ acting on $(\J,\Om,\mathbb{J})$ and study the resulting object to obtain the main results of this paper. Corollary \ref{Ham general principle} and Theorem \ref{Theorem local type} are of special note.

First, it is rather easy to prove that the space $(\J,\Om,\K)$ for the system $\K$ of paths consisting of the fiberwise geodesics is a Domic-Toledo space. In more detail for every two almost complex structures $J_0, J_1 \in \J= \Gamma(\cS; M)$ we define $[J_0,J_1]$ to be the fiberwise geodesic path $[J_0,J_1](t)$ that restricts in each fiber $\cS_x$ over a point $x \in M$ to the unique geodesic $[(J_0)_x,(J_1)_x](t)$ in $(\cS_x,\sigma_x,j_x)$ joining $(J_0)_x$ and $(J_1)_x$. Moreover, for any three elements $J_0,J_1,J_2 \in \J$ we choose $\Delta=\Delta(J_0,J_1,J_2)$ to be the fiberwise geodesic convex hull of $J_0,J_1,J_2$ so that in each fiber $\cS_x$ over $x \in M$, $\Delta$ restricts to a geodesic $2$-simplex $\Delta_x$ with respect to $\sigma_x$ with vertices $(J_0)_x,(J_1)_x,(J_2)_x$. Then since \[\int_{\Delta{(J_0,J_1,J_2)}} \Om = \int_M (\int_{\Delta_x} \sigma_x)\; \om^n(x),\] we estimate \[|\int_{\Delta{(J_0,J_1,J_2)}} \Om | \leq \int_M |\int_{\Delta_x} \sigma_x|\; \om^n(x) \leq Vol(M,\om^n) C_{\cS_n},\] as $(\cS_x,\sigma_x,j_x)$ is a Domic-Toledo space (with geodesics for the system of paths) with the constant $C_{\cS_n}$. And surely, $\J$ is contractible so the the conditions $\pi_1(J)=0$ and $\pi_2(J)=0$ are satisfied.

Second, we show that $\G = Ham(M,\om)$ is Hamiltonian-Hermitian with its action on $(\J,\Om,\K)$. First, as explained above it acts on $\J$ preserving $\Om$ with an equivariant moment map. It is also easy to deduce from the fact that the action preserves $\mathbb{J}$ that it also preserves $\K$ - though we give a direct proof. Indeed this follows immediately from the fact that for every diffeomorphism $f \in \G$ and for all $x \in M$ the map $\cS_{f^{-1} x} \to \cS_{x}$ given by $J_{f^{-1} x} \mapsto ({f_*}_x) J_{f^{-1} x} ({f_*}_x)^{-1}$ is an isometry of the Siegel upper half-spaces. The canonical metric $\rho_y$ on $\cS_{y}$ for $y \in M$ is given by $(\rho_y)_{J_y}(A_y,B_y) = const \cdot \text{trace}(A_y B_y)$ for $A_y, B_y \in T_{J_y}\cS_{y}$ ($J_y \in \cS_y$), and surely, trace is preserved by congugation with a linear isormophism.

Therefore by Theorem \ref{General principle for Hamiltonian-Hermitian} the group $\til{\G}$ admits a homogenous quasimorphism, which we show to be non-trivial by computing its local type in Theorem \ref{Theorem local type}.

\begin{cor}\label{Ham general principle}
The universal cover $\til{\G}$ of the group of Hamiltonian diffeomorphisms $\G = Ham(M,\om)$ of an arbitrary closed symplectic manifold $(M,\om)$ admits a non-trivial homogenous quasimorphism $\fS:\til{\G} \to \R$.
\end{cor}

 By construction the restriction $\fS|_{\pi_1(\G)}$ equals $\A_\mu:\pi_1(Ham(M,\omega))\to \R$. In more detail, for an element $\phi=[\{\phi_t\}] \in \pi_1(Ham(M,\omega))$ with mean-normalized Hamiltonian $H_t \in C^\infty_0(M,\R),$ we have $$A_\mu(\phi) = \int_{D}\Om - \int_0^1 dt \int_M S((\phi_t)_*J)H_t(x)\om^n,$$ where $J \in \J$ is an arbitrary element and $D$ is a disk in $\J$ spanning the loop $\{(\phi_t)_* J\}_{t \in \R/\Z}$. We now compute the homomorphism $\A_\mu$ in terms of a previously known homomorphism on $\pi_1(Ham(M,\om))$ \cite{LMP}.

\begin{df}(The homomophism $I_{c_1}:\pi_1(Ham(M,\om)) \to \R$)\label{generalized Action-Maslov c_1}
As usual with topological groups, there is a bijective correspondence between $\pi_1(Ham(M,\om))$ and the isomorphism classes of bundles $P \xrightarrow{M} S^2$ over the $2$-sphere with fiber $M$, such that their structure group is contained in $Ham(M,\om)$ \cite{LMP}. Such bundles are called \emph{Hamiltonian fiber bundles (or fibrations)} over the $2$-sphere. Over such a bundle, the vertical tangent bundle $T_V P$ is naturally endowed with the structure of a symplectic vector bundle. Hence it has Chern classes, called the \emph{vertical Chern classes}, of which we shall use the first $c^V_1:=c_1(T_V P)$. There is also a natural characteristic class $u \in H^2(P,\R)$ of such bundles with the defining properties $u|_{fiber} = [\om]$ and $\int_{fiber} u^{n+1} = 0$ (or in the case when the base is $2$-dimensional $u^{n+1}=0$) - cf. \cite{LMP,PolterovichLoops} and references therein. It is called the \emph{coupling class} of the Hamiltonian fibration. With these two characteristic classes we compose the monomial $c_1^V u^n$, where $n = \frac{1}{2} dim M$ and integrate over $P$. This yields a homomorphism $\pi_1(Ham(M,\om)) \to \R$ that we denote $I_{c_1}$. The formula for $I_{c_1}(\gamma)$ for a loop $\gamma$ in $Ham(M,\om)$ based at $Id$ is therefore \[ I_{c_1}(\gamma) = \int_{P_\gamma} c_1^V u^n, \] where $P_\gamma$ is the Hamiltonian fibration corresponding to $\gamma$.
\end{df}

\begin{thm}\label{Theorem Equality of A_mu and I_c1}
The two homomorphisms $A_\mu$ and $I_{c_1}$ from $\pi_1(\G)$ to the reals are equal.
\end{thm}

\begin{rmk}
Assume now that the almost complex structure $J_0$ is integrable - that is $(M,\om,J_0)$ is a Kahler manifold. Note that the restriction $\iota^*A_\mu$ of $A_\mu$ to the $\pi_1$ of the finite dimensional compact Lie subgroup $K:=\G_{J_0}$ of $\G$ consisting of Hamiltonian biholomorphisms satisfies $\iota^*A_\mu = -F$, for the Futaki invariant $F$ \cite{f} since the filling disk $D$ can be chosen to be trivial. The equality is understood via the isomorphism $\pi_1(K) \otimes_\Z \R \cong Lie(K)/[Lie(K),Lie(K)]$ which holds by a classical result of Chevalley and Eilenberg \cite{EilenbergChevalley} (a short account can be found in \cite{BottGr}). The consequence of Theorem \ref{Theorem Equality of A_mu and I_c1} that $I_{c_1}$ restricts to the (Bando-)Futaki invariant on $\G_{J_0}$ has previously been shown in \cite{LoopRemarks} using methods of equivariant characteristic classes.
\end{rmk}

As a corollary we answer a question of Polterovich (cf. \cite{LoopRemarks}, Discussion and Questions, 2).

\begin{cor}\label{corollary extension of I_c1 on pi1 to the quasimorphism}
We have the equality $\fS|_{\pi_1(\G)} \equiv I_{c_1}$ on $\pi_1(Ham(M,\om))$.
\end{cor}

By Proposition \ref{general symp independence} we have that $\fS$ is $Symp(M,\om)$-invariant.

\begin{cor}\label{symp invariance}
The quasimorphism $\fS: \til{\G} \to \R$ is invariant with respect to conjugation by elements of $\til{Symp}(M,\om)$ or equivalently by elements of $Symp(M,\om)$.
\end{cor}

Moreover we compute the local type of the quasimorphism $\fS$. To state the result of our computation we would first like to make two definitions of the more classical invariants in terms of which we express the answer.

\begin{df}(Calabi homomorphism on $\G_B=Ham_c(B^{2n},\om_B)$ \cite{CalabiHomomorphism}, cf. \cite{IntroSymp,PyThesis})
Given a Hamiltonian isotopy $\paph{\phi_t} \subset Ham_c(B^{2n},\om_B)$ starting at $\phi_0 = Id$ with endpoint $\phi=\phi_1$ with generating path of vector fields $\paph{X_t}$, define $H_t$ (for each $t \in [0,1]$) to be the function that vanishes near $\del B$ and satisfies $i_{X_t}\om = - d H_t.$ Then the Calabi homomorphism is defined as \[Cal_B(\paph{\phi_t}) = \intoi \int_B H_t \om^n dt.\]
It is, as can be verified using the differential homotopy and the cocycle formulas, a well-defined homomorphism $\til{\G}_B \to \R$. Moreover it vanishes on loops in $\G_B$ hence descending from $\til{\G}_B$ to $\G_B$ itself.
\end{df}

\begin{rmk}
We present a short proof that $Cal_B$ vanishes on loops in $\G_B$ that differs slightly from the one usually found in the literature. It is well-known cf. \cite{IntroSymp,PyThesis} that the Calabi homomorphism can be reinterpreted as \[Cal_B(\paph{\phi_t}) = -\frac{1}{n} \cdot \intoi \int_B (i_{X_t}\lambda)\; \om^n dt,\] for a primitive $\lambda$ of $\om_B$ in $B$. Hence \[Cal_B(\paph{\phi_t}) = -\frac{1}{n+1} \cdot \intoi \int_B (i_{X_t}\lambda - H_t)\; \om^n dt.\] However for a loop $\paph{\phi_t}$ this is proportional to \[\int_B (\int_{\paph{\phi_t x}}\lambda - \intoi H_t(\phi_t x) dt)\; \om^n(x)\] wherein the integrand is independent of $x$, as it is the Hamiltonian Action of the periodic orbit $\paph{\phi_t x}$ of $\paph{\phi_t}$. Consequently the integral localizes (up to a multiplicative constant) to the value of the integrand at each point $x \in B$ that vanishes \[\int_{\paph{\phi_t x}}\lambda - \intoi H_t(\phi_t x) dt = 0 \] for $x$ close enough to $\partial B$.
\end{rmk}

\begin{df}(The Barge-Ghys average Maslov quasimorphism on $\G_B=Ham_c(B^{2n},\om_B)$ \cite{BargeGhysEulerMaslov})\label{Barge-Ghys average Maslov quasimorphism}
Given a Hamiltonian isotopy $\paph{\phi_t} \subset Ham_c(B^{2n},\om_B)$ starting at $\phi_0 = Id$, choosing a trivialization $\Theta$ of the tangent bundle $(TB,\om_B) \cong B \times (V,\om_0)$ over $B$ as a symplectic vector bundle (here $(V,\om_0)$ is a certain symplectic vector space e.g. $(T_b B, (\om_B)_b))$ for some $b \in B$), we obtain from the family of paths of differentials $\paph{{{\phi_t}_{*_x}}:T_x B \to T_{\phi_t x} B}$ (as $x$ ranges over $B$) a family $\paph{A(x,t) \in Sp(V,\om_0)}$ of paths of symplectic linear automorphisms of $(V,\om_0)$. For each $x \in B$ we compute the value $\tau_{Lin}(\paph{A(x,t)})$ on the path $\paph{A(x,t)}$ of the Maslov quasimorphism on the universal cover of the symplectic linear group. Then the map \[\tau_{\Theta,B}: \paph{\phi_t} \mapsto \int_B \tau_{Lin}(\paph{A(x,t)}) (\om_B)^n(x)\] does not depend upon homotopies of $\paph{\phi_t}$ with fixed endpoints and yields a quasimorphism $\tau_{\Theta,B}:\til{\G}_B \to \R$. The Barge-Ghys average Maslov quasimorphism $\tau_B : \til{G} \to \R$ is its homogeneization \[\tau_B(\til{\phi}) = \displaystyle \lim_{k \to \infty} \frac{1}{k}\tau_{\Theta,B}(\til{\phi}^k).\] It does not depend on the choice of the symplectic trivialization $\Theta$. Both $\tau_{\Theta,B}$ and $\tau_B$ vanish on loops in $\G_B$ and therefore descend to quasimorphisms $\G_B \to \R$.
\end{df}

\begin{rmk}
The vanishing of $\tau_{\Theta,B}$ on loops can be shown by a similar localization argument as for the Calabi homomorphism. Indeed for a loop $\paph{\phi_t}$ in $\G_B$ the value $\tau_{Lin}(\paph{A(x,t)})$ equals the Maslov index of the loop $\paph{A(x,t)}$ which by the homotopy invariance
of the Maslov index is independent of $x$, and for $x$ near $\partial B$ the loop $\paph{A(x,t)}$ is trivial. Hence the integrand vanishes for all $x \in B$, wherefrom $\tau_{\Theta,B}(\paph{\phi_t}) = 0$.
\end{rmk}

\begin{thm}\label{Theorem local type}
 Let $c = n \int_M c_1 \om^{n-1}/ \int_M \om^n = \int S(J) \om^n / Vol(M,\om^n)$ be the average Hermitian scalar curvature. Then the restriction of $\fS$ to the subgroup $\G_B = Ham_c(B,\om|_B) \subset \G$ of Hamiltonian diffeomorphisms supported in an embedded ball $B$ in $M$ satisfies \[\fS|_{\G_B} = \frac{1}{2}\tau_B - c \text{Cal}_B,\] where $\tau_B$ is the Barge-Ghys Maslov quasimorphism on $\G_B = Ham_c(B^{2n},\om_{std})$ and $\text{Cal}_B$ is the Calabi homomorphism.
\end{thm}

We describe the relation of the quasimorphism $\fS$ to the quasimorphisms $\fS_{Py}$ and $\fS_{En}$ introduced by Py \cite{pyqm,PyThesis} for closed manifolds $(M,\om)$ with $c_1(TM,\om) = \kappa [\om]$ for $\kappa \neq 0$ and Entov \cite{EntovCommutatorLength} for closed manifolds $(M,\om)$ with $c_1(TM,\om)=0$. First we state briefly the definitions of the quasimorphisms $\fS_{Py}$ and $\fS_{En}$. The detailed definitions appear in the proofs section.

\begin{df}\label{Py quasimorphism}(A sketch of a definition of $\fS_{Py}$ \cite{pyqm,PyThesis})
Endow the unit frame bundle $P \xrightarrow{S^1} M$ of $L=\Lambda^n_\C(TM,J,\om)$ for a compatible complex structure $J \in \J$ with the structure $\alpha_0$ of a prequantization of $(M,-\om)$. Note that there is a natural map $det^{2}: \cL(TM,\om) \to P^2$ from the Lagrangian Grassmannian bundle $\cL(TM,\om)$ to the unitary frame bundle $P^2$ of $L^{\otimes 2}$, since $\cL(TM)_x = U(TM_x,\om_x,J_x)/O(n)$. Note that $\alpha_0$ induces a structure $\alpha$ of prequantization of $(M,-2\om)$ on $P^2$. Given a path $\arr{\phi} = \paph{\phi_t}$ in $\G$ with $\phi_0 = Id$, choosing a point $L \in \mathcal{L}(TM,\omega)_x$ we have the curve $\{{{\phi_t}_*}_x(L)\}_{0 \leq t \leq 1}$ in $\mathcal{L}(TM,\omega)$ and considering $\arr{\phi}$ as a path of Hamiltonian isotopies of $(M,-2\omega)$ we have the canonical lifting $\{\widehat{\phi}_t\}_{0 \leq t \leq 1},\;\widehat{\phi}_0 = Id$ of $\arr{\phi}$ to the identity component $Q=Quant(P^2,\alpha)$ of the group of diffeomorphisms of $P^2$ that preserve $\alpha$. Consequently, one considers the two curves \[\{ det^2({\phi_t}_{*_x}(L))\}_{0 \leq t \leq 1}\]  and  \[\{\widehat{\phi}_t(det^2(L))\}_{0 \leq t \leq 1}\] in $P^2$. Both these curves in $P^2$ start at $det^2(L)$ and cover the path $\{\phi_t (x)\}_{0 \leq t \leq 1}$ in $M$ and hence differ by an angle: $$det^2({\phi_t}_*(L)) = e^{i2\pi \vartheta(t)} \widehat{\phi}_t(det^2(L)),$$ for a continuous function $\vartheta:[0,1] \to \R$. Define a continuous function on $\mathcal{L}(TM,\omega)$ by $$angle(L,\arr{\phi}):=\vartheta(1) - \vartheta(0).$$ Then the function $angle(x,\arr{\phi}) =  \displaystyle\inf_{L \in \mathcal{L}(TM,\omega)_x} angle(L,\arr{\phi})$ on $M$ is measurable, bounded and defines the quasimorphism \[S_2(\arr{\phi}) = -\int_M angle(x,\arr{\phi}) \omega^n(x)\] that does not depend upon homotopies of $\arr{\phi}$ with fixed endpoints and is thus defined as a real-valued function on $\til{G}$.
Its homogeneization $\mathfrak{S}_{Py}:\til{G} \to \R$, defined by $\mathfrak{S}_{Py}(\til{\phi}):=\displaystyle\lim_{k \to \infty} \frac{S_2(\til{\phi}^k)}{k}$ is a homogenous quasimorhism on $\til{G}$ that is independent of the non-canonical structure on $P$ of a prequantization of $(M,-\omega)$, of the prequantization form $\alpha$ on it and of the almost complex structure $J$.
\end{df}

\begin{df}\label{Entov quasimorphism}(A sketch of a definition of $\fS_{En}$ \cite{EntovCommutatorLength})
Given a symplectic manifold $(M,\om)$ with $c_1(TM,\om) = 0$ one first trivializes $(TM,\om,J)$ for $J \in \J$ as a Hermitian vector bundle over the complement $U = M\setminus Z$ of a compact triangulated subset $Z$ of $\codim(Z) \geq 3$, where the differential of the trivialization, appropriately defined, is uniformly bounded. For a path $\arr{\phi}=\paph{\phi_t}$ in $\G$ by relaxing $Z$ to be a countable union $Z_{\arr{\phi}} = \bigcup_{j \in \Z} Z_j$ depending on $\arr{\phi}$ of sets $Z_j$ of $\codim(Z_j) \geq 2$ one can assume that $U$ is invariant with respect to $\phi_t$ for all $t$. Then from the path $\paph{{\phi_t}_{*_x}}$ for $x \in U$ one obtains a continuous path $\paph{A(x,t)}$ with $A(x,0) = Id$ in $Sp(2n,\R)$ and proceeds to define \[angle(x,\arr{\phi}) = varangle({\paph{det^2(A(x,t))}}).\] One then shows that this function extended by $0$ on $Z$ is integrable on $M$ and that \[T_1(\arr{\phi}) = \int_M angle(x,\arr{\phi}) \om^n(x)\] does not depend on homotopies of $\arr{\phi}$ with fixed endpoints (by relaxing $Z$ to be a countable union of sets $Z'_j$ of $\codim(Z'_j) \geq 1$ depending on a given homotopy), and defines a quasimorphism \[T_1: \til{\G} \to \R.\]
Its homogeneization $\mathfrak{S}_{En}:\til{G} \to \R$, defined by $\mathfrak{S}_{En}(\til{\phi}):=\displaystyle\lim_{k \to \infty} \frac{T_1(\til{\phi}^k)}{k}$ is a homogenous quasimorhism on $\til{G}$ that is independent of the non-canonical choices of trivialization, of the set $U = M\setminus Z$ and of the almost complex structure $J$.
\end{df}

We claim that the quasimorphism $\til{G} \to \R$ obtained from Corollary \ref{Ham general principle} agrees with these two quasimorphisms in the settings of their definitions.

\begin{thm}\label{General agrees with Py and Entov}
1. On symplectic manifolds $(M,\om)$ with $c_1(TM,\om) = \kappa [\om]$ for $\kappa \neq 0$ we have $2\fS=\nobreak-\fS_{Py}.$
2. On symplectic manifolds $(M,\om)$ with $c_1(TM,\om) = 0$ we have $2\fS=\nobreak\fS_{En}.$
\end{thm}

\begin{rmk}
The analogues of Theorem \ref{Theorem local type} for the cases 1. and 2. above were shown in \cite{EntovCommutatorLength,pyqm}. The analogue of Corollary \ref{corollary extension of I_c1 on pi1 to the quasimorphism} was shown in \cite{PyThesis}. The agreement of our results with the ones shown in these papers is as follows. For analogues of Theorem \ref{Theorem local type} note that the average scalar curvature $c$ satisfies $c = n\kappa$ when $c_1(TM,\om) = \kappa [\om]$, for every $\kappa$. For the analogue of Corollary \ref{corollary extension of I_c1 on pi1 to the quasimorphism}, use the easy Computation 1 from \cite{LoopRemarks}, near the end of Section 1.2.
\end{rmk}

\begin{rmk}
We would also like to note that the general scheme of Theorem \ref{General principle for Hamiltonian-Hermitian} applies to the construction of quasimorphisms $\til{\G} \to \R$ for the group $\G = Ham_c(M,\om)$ of Hamiltonian diffeomorphisms with compact support of symplectic manifolds $(M,\om)$ of finite volume (without boundary) that are not compact. Indeed, $\J$ here is also a Domic-Toledo space, since $\int_M \om^n$ is finite, and Donaldson's theory for the scalar curvature as an equivariant moment map \cite{DonaldsonGaugeTheoryComplexGeometry} applies here nearly verbatim. The only difference is that the symplectic form $\Om$ is not defined on all the tangent space $T_{J_0} \J$ - indeed given $A,B \in T_{J_0} \J$ the function $\sigma_{(J_0)_x}(A_x,B_x)$ may well be non-integrable with respect to $\om^n$. However, since we compute for diffeomorphisms with compact support, all relevant computations happen in a compact subset of $M$ where all functions that appear are integrable. Moreover, all functions, vector fields, one-forms and sections of endomorphism bundles have compact support, therefore the only non-local part in Donaldson's proof \cite{DonaldsonGaugeTheoryComplexGeometry} - integration by parts to show the actual integral formulae - goes through (all the other arguments are local). At the same time, when the symplectic volume of $M$ is not finite, $\J$ stops being a Domic-Toledo space (at least with the natural definitions) and hence this approach does not seem to give quasimorphisms. It would be interesting to investigate the restriction to $\pi_1(\G)$ of the quasimorphism in the finite volume case. The local type is obtained by nearly the same computation as the one given for the closed case and is given by the Barge-Ghys average Maslov quasimorphism $\tau$. The $\cH=Symp_c(M,\om)$-invariance holds as before.
\end{rmk}

A corollary, as obtained in \cite{EntovCommutatorLength} for symplectic manifolds with $c_1(TM,\om) = 0,$ is that the commutator length of $\til{\G}$ is unbounded.

\begin{cor}
The diameter in the commutator length of the group $\til{\G}$ for $\G = Ham(M,\om)$ of a closed symplectic manifold $(M,\om)$ and of the perfect (\cite{BanyagaPerfectSimple}) group $\mathrm{Ker}(Cal: \til{G} \to \R)$ for $\G=Ham_c(M,\om)$ of an open finite volume symplectic manifold $(M,\om)$ is infinite. In the closed case, under the additional assumption $I_{c_1} \equiv 0$, the same conclusion follows for $\G$ itself.
\end{cor}

We also note that for reasons of naturality of the constructions and normalizations of Hamiltonians we have the following proposition.

\begin{prop}\label{embedding functoriality}(Embedding functoriality)
Given an open subset $U \subset M$ of a closed symplectic manifold $(M,\om)$, denote by $\mathfrak{S}_M$ the quasimorphism obtained on $\G$ for $(M,\om)$ and by $\mathfrak{S}_U$ the quasimorphism obtained on $\G_U$ for $(U,\om|_U)$. Then $$\mathfrak{S}_M|_{\G_U} = \mathfrak{S}_U - c \cdot Cal_U,$$ for the average Hermitian scalar curvature $c$. Similarly if $M$ were an open symplectic manifold of finite volume $$\mathfrak{S}_M|_{\G_U} = \mathfrak{S}_U.$$
\end{prop}

\subsection{Application to the $L^2_2$-distance on $\til{Ham}(M,\om)$}

While it is not surprising that our quasimorphism is bounded by a multiple of the Sobolev $L^2_2$ norm on $\til{Ham}(M,\om)$, indeed $\fS_{J_0}$ is surely continuous in the $C^1$-topology induced to $Ham(M,\om)$ from $Diff(M)$, we present a proof for the sheer simplicity of the argument.

For a Hamiltonian isotopy $\arr{\phi}=\paph{\phi_t}$ of a symplectic manifold $(M,\om)$ starting at the identity that is generated by the zero-mean-normalized Hamiltonian $H_t$ put \[||\arr{\phi}||_{k,p} = \intoi ||H_t||_{L^p_k(M,\om^n)}.\] Then define the norm of an element $\til{\phi} \in \til{\G}$ by \[||\til{\phi}||_{k,p} = \displaystyle\inf_{[\arr{\phi}] = \til{\phi}} ||\arr{\phi}||_{k,p}.\] Finally define the norm of $\phi \in \G$ by $||\phi||_{k,p} = \displaystyle\inf_{\pi(\til{\phi}) = \phi} ||\til{\phi}||_{k,p}$ for the natural projection $\pi: \til{\G} \to \G$. For two elements $a,b$ of the above groups define the distance \[d_{p,k}(a,b) = ||a^{-1}b||_{k,p}.\]  The following facts are easy to check.
\begin{itemize}
\item For $k \geq 1$ the $(p,k)$-norms and distances are equivalent to $(p,k-1)$-norms and distances as defined via the vector field $X_t$ generating $\arr{\phi}$.
\item For $k \geq 1$ these norms and distances are non-degenerate.
\end{itemize}

We show in Section \ref{calibrating the L^2_2 norm} that $\fS_{J_0}$ calibrates the $(2,2)$-norm as follows

\begin{equation}\label{bound by L^2_2 of fS} \fS_{J_0}(\til{\phi}) \leq C(n,\om,J_0) ||\til{\phi}||_{2,2},\end{equation}

for a constant $C(n,\om,J_0)$ that does not depend on $\til{\phi}$. As a corollary we obtain that the $L^2_2$ distance is unbounded on $\til{\G}$.

\begin{cor}
The diameter of $\til{\G}$ is infinite with respect to the $L^2_2$-distance for every symplectic manifold $(M,\om)$ of finite volume.
\end{cor}

\begin{rmk}
Given that $\A_\mu: \pi_1 \G \to \R$ vanishes, the same consequence holds for the $L^2_2$-distance on the group $\G$ itself. For closed manifolds this condition is equivalent to the vanishing of $I_{c_1}$.
\end{rmk}

 The unboundedness of the $L^2_1$-metric on compact exact symplectic manifolds was previously proven by Eliashberg and Ratiu \cite{EliashbergRatiu} (their methods work even for the larger group $\cH = Symp(M,\om)$ with appropriate definitions), while sharper topological bounds for the $2$-disc were obtained by Gambaudo and Lagrange \cite{GambaudoLagrange} (cf. \cite{BenaimGambaudoMetricProperties}, \cite{BrandenburskyEstimates}).

\subsection{Finite dimensional examples: Guichardet-Wigner quasimorphisms}\label{FinDimExamples}

The general principle outlined in Section \ref{quasimorphisms} applies also to finite dimensional Hermitian Lie groups acting on their corresponding Hermitian symmetric spaces of non-compact type. In this section we describe this application, in part for use in the proofs later.

Let $G$ be a simple Hermitian symmetric Lie group. Then the adjoint form of $G$ belongs to those of the following list of Lie groups: $SU(p,q),\,SO_0(2,q)\;q\neq 2,\, Sp(2n,\R),\, SO^*(2n)\;n\geq2$ and two real forms of the complex simple Lie groups of types $E_6$ and $E_7$ respectively. Let us assume that the center of $G$ is finite, so that $\pi_1(G)$ is infinite. Let $K \subset G$ be the analytic subgroup corresponding to the maximal compact Lie subalgebra $\mathfrak{k}$ of $\mathfrak{g}$. In this situation there is a corresponding Hermitian symmetric space $X = G/K$, endowed with a natural complex structure $j_X$ and a Kahler form $\sigma_X$ that is invariant with respect to the transitive action of the group $G$ (proportional to the Bergman Kahler structure when such a space is realized as a symmetric bounded domain in a complex affine space by the Harish-Chandra embedding) cf. \cite{Helgason,Knapp,Mok}. The works of Domic-Toledo and {\O}rsted \cite{DomicToledo,ClercOrsted} show that when we take the system of paths $\K$ to consist of the geodesics with respect to the invariant Kahler metric, then $(X,\sigma_X,\K)$ is a Domic-Toledo space in our terminology.

Moreover, we note that by e.g. \cite{Mok} these spaces $(X,\sigma_X,j_X)$ are Kahler-Einstein manifolds (that is - their Ricci forms are proportional to their Kahler forms: $Ric(\sigma_X) = \lambda \sigma_X$, where for the Bergman metric we have $\lambda = -1$). Note that $Ric(\sigma_X)$ is equal up to a universal constant to the curvature of the Chern connection on the line bundle $L_X=\Lambda_\C^N TX$, with the holomorphic and Hermitian structures induced by $j_X$ and $\sigma_X$.

We now show that $G$ is Hamiltonian-Hermitian with its action on $X$. Firstly the group $G$ acts on $X$ by maps preserving $j_X$ and $\sigma_X$ (symplectic biholomorphisms) and hence preserving the system of geodesics $\K$. We now claim that the group $G$ acts on $X$ with an equivariant moment map $\mu_X: Lie(G) \times X \to \R$. Note that as the Chern connection on $TX$ is given canonically by $(\sigma_X, j_X)$ and the action preserves these structures, it will also preserve the Chern connection. Consider the natural lift of the action of $G$ on $X$ to an action of $G$ on $TX$ by taking differentials. This induces an action of $G$ on $L_X=\Lambda_\C^N TX$. Note that this action preserves the Hermitian structure on $L_X$, and hence it descends to the circle bundle $P_X \xrightarrow{S^1} X$ of unit vectors in $L_X$ (the unitary frame bundle of the Hermitian vector bundle $L_X$). The Chern connection on $L_X$ induces a real-valued connection one-form (cf. \cite{LoopRemarks} Appendix A) $\alpha_X$ on the principal $S^1$-bundle $P_X$ over $X$, that by the Kahler-Einstein property satisfies the relation \begin{equation}\label{Chern prequantization for Hermitian symmetric spaces of non-compact type}d \alpha_X = \til{\sigma_X}\end{equation} for the lift $\til{\sigma_X}$ of $\sigma_X$ to $P_X$ by the natural projection $P_X \to X$, as follows from what is noted above. Now the action of the group $G$ on $P_X$ covering the action of $G$ on $X$ preserves the one-form $\alpha_X$ (by preservation of the Chern connection). This is enough to give an equivariant moment map for the action of $G$ on $X$. Indeed, it is constructed as follows. A vector $\xi \in Lie(G)$ induces the vector field $\bar{\xi}$ on $X$ by the action of $G$ on $X$ and a vector field $\widehat{\xi}$ on $P_X$ that covers $\bar{\xi}$, by the action of $G$ on $P_X$. We claim that the equivariant moment map is given by \begin{equation}\label{equation - moment map for finite dimensional groups}\mu_X(\xi)(x) = (\alpha_X)_y(\widehat{\xi}_y)\end{equation} for any $y \in P_X$ over $x \in X$ (indeed $\widehat{\xi}$ is equivariant with respect to the natural circle action on $P_X$ as is $\alpha_X$ and hence $(\alpha_X)_y(\widehat{\xi}_y)$ does not depend on the choice of $y$ over $x$). Firstly by relation (\ref{Chern prequantization for Hermitian symmetric spaces of non-compact type}) and the preservation $\cL_{\widehat{\xi}} \alpha_X = 0$ of the connection by the infinitesimal action we have \[i_{\bar{\xi}} \sigma_X = - d \mu(\xi)(x).\] Hence $\mu_X$ is a moment map for the action of $G$ on $X$. For the equivariance we note once again that the action of $G$ on $P_X$ preserves $\alpha_X$ and that the vector field $\widehat{\xi}$ has a corresponding equivariance property. Namely for any $g \in G$ and $y \in P$ denoting by $\wh{g} \cdot y$ the action of $G$ on $P_X$ and by $\wh{g}_{*_y}$ the corresponding differential $T_y P_X \to T_{\wh{g}\cdot y}P_X$ we have the very general equivariance property for infinitesimal actions corresponding to Lie group actions on spaces \[\wh{\xi}(\wh{g} \cdot y) = \wh{g}_{*_y} (\wh{Ad_{g^{-1}}\xi}(y)).\] Now noting that for $y \in P_X$ over $x \in X$ the point $\wh{g} \cdot y$ is over $g \cdot x$, we obtain \[\mu_X(Ad_g \xi)(g\cdot x) = (\alpha_X)_{\wh{g}\cdot y}(\widehat{Ad_{g} \xi}({\wh{g} \cdot y})) = (\alpha_X)_{\wh{g}\cdot y}( \wh{g}_{*_y} \wh{\xi}(y)) = (\alpha_X)_{y}(\wh{\xi}(y)) = \mu_X(\xi)(x),\] showing equivariance.

Hence $G$ is Hamiltonian-Hermitian with Domic-Toledo space $(X,\sigma_X,j_X)$ and equivariant moment map $\mu_X$, and therefore by Theorem \ref{General principle for Hamiltonian-Hermitian} has a homogenous quasimorphism.

\begin{cor}\label{finite dimensional nu}
Theorem \ref{General principle for Hamiltonian-Hermitian} gives a homogenous quasimorphism $\nu_G: \til{G} \to \R$ for every simple Hermitian symmetric Lie group $G$.
\end{cor}

It remains to show that it is non-trivial. In fact we show in Section \ref{finite dimensional and Guichardet-Wigner} that it is equal to the Guichardet-Wigner \cite{GuichardetWigner,DupontGuichardet,ClercKoufany,ShternAutomaticContinuity,BurgerIozziWienhardMaximalToledoInvariant} quasimorphism $\varrho_G$ on $\til{G}$ by comparing them on $\pi_1(G)$ and arguing that a homogenous quasimorphism on $\til{G}$ is determined by its restriction to the fundamental group.

\begin{prop}\label{finite dimensional nu and varrho GW}
The quasimorphisms $\nu_G$ and $\varrho_{G}$ on $\til{G}$ satisfy the equality \[\nu_G = -\varrho_{G}.\]
\end{prop}

We would now like to give a reformulation of the construction of $\nu_x$ in the finite dimensional case as a certain rotation number. Indeed consider once again the principal $S^1$-bundle $P_X \to X$. Trivialize it by taking parallel transports $\Gamma_{\gamma(y,x)}: (P_X)_y \to (P_X)_x$ along geodesics $\gamma(y,x)$ for $y \in X$. Then given a path $\arr{g}=\paph{g_t}$ with $g_0 = Id$ in $G$, the path of differentials $(g_t)_{*_x}:T_x X \to T_{g_t\cdot x}X$ gives us a path $\Gamma_{\gamma(g_t\cdot x,x)} \circ \widehat{g}_t|_{(P_X)_x}: (P_X)_x \to (P_X)_x$ which we consider as a path in $U(1) \cong S^1$. Then \begin{equation}\label{nu_x via varangles}\nu_x(\arr{g}) = \text{varangle}(\paph{\Gamma_{\gamma(g_t\cdot x,x)} \circ \widehat{g}_t|_{(P_X)_x}}).\end{equation} Indeed, denoting $\gamma_t:=\gamma(g_t\cdot x,x)$ and $\beta_t = \{g_{t'}\cdot x\}_{t'=0}^t$ we have \[varangle(\paph{\Gamma_{\gamma(g_t\cdot x,x)} \circ \widehat{g}_t|_{(P_X)_x}}) = varangle(\paph{\Gamma_{\gamma_t} \circ \Gamma_{\beta_t}}) + varangle(\paph{\Gamma_{\overline{\beta}_t} \circ \widehat{g}_t|_{(P_X)_x}}) = \] \[ = \int_{D_{\arr{g}}} \sigma_X - \intoi (\alpha_X)_{{\wh{g}}_t\cdot y}(\widehat{\xi}_t)_{{\wh{g}_t}\cdot y} dt = \int_{D_{\arr{g}}} \sigma_X - \intoi \mu(\xi_t)(g_t\cdot x) dt =\nu_x(\arr{g}).\]

It is interesting to note that taking this reformulation as a definition for the quasimorphism, its independence upon homotopies with fixed endpoints follows immediately by continuity.

\subsection*{Acknowledgements}
First and foremost I thank my advisor Leonid Polterovich for his support and encouragement, for his continuous interest in this project and for many fruitful discussions. I have also benefited from several stimulating conversations with Pierre Py. A decisive part of this project was carried out during the author's visit to the Mathematics Department at the University of Chicago. I thank Leonid Polterovich and the Mathematics Department for their hospitality and for a great research atmosphere. I thank Akira Fujiki for sending me a proof of his theorem (Equation \ref{FujikiFormula}). I also thank Marc Burger and Danny Calegari for useful comments. Many thanks are due to the referee for remarks that have improved the exposition. This paper is partially supported by the Israel Science Foundation grant $\#$\nobreak509/07.

\section{Proofs}\label{Proofs}

\subsection{The Action homomorphism}
We prove that the number $\A_\mu(a)$ defined in Section \ref{action} is well-defined and determines a homomorphism $\pi_1(\G) \to \R/\Pe_\Om$. We refer to Sections \ref{section - Moment maps} and \ref{action} for the relevant notation and definitions. Let us first prove that it is well-defined. First of all, the value $\A_\mu(a) \in \R/\Pe_\Om$ obviously does not depend on the spanning disk.
Let us prove that it does not depend on the point $x \in \X$. Take another point $x' \in \X$ and choose a path $\beta : [0,1] \to \X$ between the two: $\beta(0) = x, \; \beta(1) = x'.$ Consider the cylindric cycle $C:\R/\Z \times [0,1] \to \X$ defined by $C(t,s) = \phi_t \cdot \beta(s)$. Note that $C(t,0) = \phi_t \cdot x$ and that $C(t,1) = \phi_t \cdot x'$. Define a spanning disk $D'$ for $\phi_{x'}$ by $D' = D \cup_{\phi_x} C$. Then the equality $$\int_D \Omega - \int_0^1\mu(X_t)(\phi_t \cdot x)\,dt = \int_{D'} \Omega - \int_0^1\mu(X_t)(\phi_t \cdot x')\,dt $$ that we are trying to prove reduces to $$ \int_{D'} \Omega - \int_D \Omega = \int_0^1\mu(X_t)(\phi_t \cdot x')\,dt - \int_0^1\mu(X_t)(\phi_t \cdot x)\,dt,$$ which is equivalent to
$$ \int_C \Omega = \int_0^1\mu(X_t)(C(t,1))\,dt - \int_0^1\mu(X_t)(C(t,0))\,dt.$$ This equality is established by direct computation of the left hand side. Indeed $$\int_C \Omega = \int_0^1 \int_0^1 \Om(\del_s C(s,t),\del_t C(s,t)) ds dt = \int_0^1 \int_0^1 \Om(\del_s C,\Xi_t(C(t,s))) ds dt =$$$$= \int_0^1 \int_0^1 d_{C(s,t)}\mu(X_t)(\del_s C) ds dt = \int_0^1 \int_0^1 \del_s \mu(X_t)(C(s,t)) ds dt =$$$$= \int_0^1 \mu(X_t)(C(1,t)) - \mu(X_t)(C(0,t)) dt = \int_0^1\mu(X_t)(C(t,1))\,dt - \int_0^1\mu(X_t)(C(t,0))\,dt,$$ yielding the desired equality.

Let us now proceed to prove that $\A_\mu(\{\phi_t\})=\int_D \Omega - \int_0^1\mu(X_t)(\overline{\phi_t} x)\,dt $ remains invariant when $\phi_t$ is deformed homotopically with fixed endpoints. Let $\phi^s_t, \; 0\leq s \leq 1$ be such a homotopy. That is $\phi^s_0 \equiv Id,\;$ $\phi^s_1 \equiv Id$ and $(s,t) \to \phi^s_t$ is a smooth map $[0,1]\times[0,1] \to \G$. Surely, it is enough to prove that for all $s$ the derivative $\frac{\del}{\del s}|_s\A_\mu(\phi^s_t)$ vanishes. To this end we use the following lemma, which is a direct consequence of the standard differential homotopy formula.

\begin{lma}
Let $X^s_t$ and $Y^s_t$ be the elements $X^s_t = \frac{\del}{\del \tau}|_{\tau = t}\phi^s_\tau \cdot (\phi^s_t)^{-1}$, $Y^s_t = \frac{\del}{\del \sigma}|_{\sigma = s}\phi^\sigma_t \cdot (\phi^s_t)^{-1}$ of $Lie (\G)$ (note that $Y^s_0 \equiv 0$ and $Y^s_1 \equiv 0$).  Then $\ddel{s} Ad_{(\phi^s_t)^{-1}} X^s_t = Ad_{(\phi^s_t)^\mone} \ddel{t} Y^s_t$ and $\ddel{t} Ad_{(\phi^s_t)^\mone} Y^s_t = Ad_{(\phi^s_t)^\mone} \ddel{s} X^s_t.$
\end{lma}

\begin{pf}
The differential homotopy formula says $\ddel{s} X^s_t = \ddel{t} Y^s_t + [X^s_t,Y^s_t]$. Differentiating $\ddel{s} Ad_{(\phi^s_t)^\mone} X^s_t$ we obtain $Ad_{(\phi^s_t)^\mone}([Y^s_t,X^s_t] + \ddel{s}X^s_t),$ which by the differential homotopy formula equals $Ad_{(\phi^s_t)^\mone} \ddel{t} Y^s_t$. The other equality is obtained in the same way. Both are equivalent to the original differential homotopy formula.
\end{pf}

Now \[-\frac{\del}{\del s}|_s\A_\mu(\phi^s_t) = \ddel{s}\int_0^1\mu(X^s_t)(\overline{\phi^s_t} x)\,dt - \ddel{s}\int_{D^s_t} \Omega =\]\[=
\int_0^1 \mu(\ddel{s} X^s_t)(\overline{\phi^s_t} x)\,dt + \int_0^1 \iota_{\Upsilon^s_t(\overline{\phi^s_t} x)} d_{\overline{\phi^s_t}x}\mu(X^s_t)\,dt - \int_0^1 \iota_{\Upsilon^s_t(\overline{\phi^s_t} x)}\Omega (\Xi^s_t(\overline{\phi^s_t} x)) =\]\[= \int_0^1 \mu(Ad_{(\phi^s_t)^\mone}\ddel{s} X^s_t)(x)\,dt + \int_0^1 \Om (\Xi^s_t(\overline{\phi^s_t} x),\Upsilon^s_t(\overline{\phi^s_t} x))\,dt - \int_0^1 \Om (\Upsilon^s_t(\overline{\phi^s_t} x),\Xi^s_t(\overline{\phi^s_t} x))=\]\[= \int_0^1 \mu(\ddel{t} Ad_{(\phi^s_t)^\mone} Y^s_t)(x)\,dt = \int_0^1 \ddel{t}\mu( Ad_{(\phi^s_t)^\mone} Y^s_t)(x)\,dt =\]\[= \mu( Ad_{(\phi^s_1)^\mone} Y^s_1)(x) - \mu( Ad_{(\phi^s_0)^\mone} Y^s_0)(x) = \mu(Y^s_1)(x) - \mu(Y^s_0)(x) = 0.\] This yields the desired equality. Here ${D^s_t}$ is obtained by gluing $D$ and $C(s,t)=\overline{\phi^s_t} x$ along $\phi_x$. The vector fields $\Xi^s_t, \Upsilon^s_t$ are the infinitesimal actions of $X^s_t$ and $Y\st$.

At last, let us prove that $\A_\mu$ defines a homomorphism $\pi_1(\G) \to \R/\Pe_\Om$. Indeed, take two loops $\phi=\{\phi_t\}$, $\psi=\{\psi_t\}$ based at $Id$. Consider their concatenation $\chi = \psi \ast \phi$. Then $\chi_x = \phi_x \ast \psi_x$. Moreover we can choose a spanning disk of $\chi_x$ that factors through the topological wedge (we refer to \cite{FomenkoFuks} for the definition of wedge and for related notations) of the spanning disks $D_\phi$, $D_\psi$ of $\phi_x$, $\psi_x$. That is $D: \D \to \X$ factors as $D: \D \xrightarrow{pr} \D \bigvee_{1} \D \xrightarrow{D_\phi \bigvee_{x} D_\psi} \X$. Hence $\A_\mu(\chi)= \int_{D_\phi} \Omega - \int_0^1 \int_0^1\mu(X_t)(\overline{\phi_t} x) + \int_{D_\psi} \Omega - \int_0^1\mu(Y_t)(\overline{\psi_t\phi_1} x)\,dt = \int_{D_\phi} \Omega - \int_0^1 \int_0^1\mu(X_t)(\overline{\phi_t} x) + \int_{D_\psi} \Omega - \int_0^1\mu(Y_t)(\overline{\psi_t} x)\,dt = \A_\mu(\phi) + \A_\mu(\psi).$ Here $X_t$ and $Y_t$ are the elements of $Lie (\G)$ corresponding to $\phi$ and $\psi$. The penultimate equality follows from the fact that $\phi_1 = Id$.

\subsection{The quasimorphism on Hamiltonian-Hermitian groups}\label{proof of general principle for Hamiltonian-Hermitian}

We now prove Theorem \ref{General principle for Hamiltonian-Hermitian} on the construction of quasimorphisms on Hamiltonian-Hermitian groups.

\begin{proof}

First, the independence on the disk follows trivially, since $\pi_2(\X) = 0$ and $\Om$ is closed.

We proceed to show that the map is independent upon homotopies of $\paph{g_t}$ with fixed endpoints. Let $g\st$ be a homotopy with fixed endpoints $g^s_0 \equiv Id$, $g^s_1 \equiv g_1$ of $\overrightarrow{g}=\paph{g^0_t = g_t}$ to $\overrightarrow{h}=\paph{g^1_t=h_t}$. Note that in this situation the concatenation $q =\paph{g_t} \#\; \overline{\paph{h_t}}$ is a contactible loop in $\G$ based at $Id$. Denote by $C$ the disk $C=\{g\st\cdot x\}_{0\leq s,t \leq 1}$. Choose the disks of integration as follows. When computing for $\paph{h_t}$ choose an arbitrary disk $D_{\paph{h_t}}$ and for $\paph{g_t}$ choose $D_{\paph{g_t}} = D_{\paph{h_t}} \cup C$ where the gluing is over the common path $\paph{h_t \cdot x}$.
Then $$\nu_{x}(\arr{g})-\nu_{x}(\arr{h}) = \int_C \Om - (\int_0^1 \mu(X_t)(g_t \cdot x) - \int_0^1 \mu(Y_t)(h_t \cdot x)) = \A_{\mu}(q) = 0$$ since $\A_\mu$ is a homomorphism on $\pi_1(\G)$ and $q$ is a contractible loop. Here $\paph{X_t}$ and $\paph{Y_t}$ are the paths in $Lie(\G)$ corresponding to $\paph{g_t}$ and $\paph{h_t}$

We now show the quasimorphism property of $\nu_{x}$. Take two paths $\arr{g}=\{g_t\}_{0 \leq t \leq 1},\arr{h}=\{h_t\}_{0 \leq t \leq 1}$ representing elements $\til{g},\til{h}$ of $\til{\G}$. Denote by $g,h$ their endpoints. We would like to compare $\nu_{x}(\til{g}\til{h})$ with $\nu_{x}(\til{g}) + \nu_{x}(\til{h})$. Note that $\til{g}\til{h}$ is represented by the path $\arr{g} \#\; g_1\arr{h}$, where $g_1 \arr{h} = \paph{g_1 h_t}$. Hence we will compare $$\nu_{x}(\arr{g} \#\; g_1\arr{h}) \;\;\;\;\;\text{to}\;\;\;\;\; \nu_{x}(\arr{g}) + \nu_{x}(\arr{h}).$$ The definition of $\nu_x$ involves two summands - one involving the symplectic area and one involving the moment map. We first show that the terms involving the moment maps are equal. And indeed \begin{equation}\label{moment map part of qm property formula}\intoi \mu(X_t)(g_t \cdot x)dt + \intoi \mu(Ad_{g_1}Y_t)(g_1 \cdot h_t \cdot x)dt = \intoi \mu(X_t)(g_t \cdot x)dt + \intoi \mu(Y_t)( h_t \cdot x)dt,\end{equation} by the equivariance of the moment map.

Now we show that the terms involving symplectic area agree up to the function \begin{equation}\label{triangle and defect}\int_{\Delta(x,g\cdot x,gh \cdot x)} \Om\end{equation} which is bounded by a constant $C_\X$ depending only on the Domic-Toledo space $(\X,\Om,\K)$. Indeed choosing arbitrary disks of integration $D_{\arr{g}}$ for $\arr{g}$ and $D_{\arr{h}}$ for $\arr{h}$, choose $\arr{g} \#\; g_1\arr{h}$ the disk $D_{\arr{g} \#\; g_1\arr{h}} = (D_{\arr{g}} \cup g_1\cdot D_{\arr{h}}) \cup \Delta(x,g\cdot x,gh \cdot x)$ where the gluing is over the common path $[x,g\cdot x]\cup g[x,h\cdot x]$, which equals $[x,g\cdot x]\cup [g\cdot x,gh\cdot x]$ by preservation of $\K$. Hence \begin{equation}\label{symplectic area part of the quasimorphism property}\int_{D_{\arr{g} \#\; g_1\arr{h}}} \Om = \int_{D_{\arr{g}}}\Om + \int_{g_1\cdot D_{\arr{h}}}\Om + \int_{\Delta(x,g\cdot x,gh \cdot x)}\Om\\ =  \int_{D_{\arr{g}}}\Om + \int_{D_{\arr{h}}}\Om + \int_{\Delta(x,g\cdot x,gh \cdot x)}\Om, \end{equation} by preservation of $\Om$ by the action. This finishes the proof of the quasimorphism property.

Now we discuss the independence of the homogeneization $$\nu(\til{g}) = \displaystyle \lim_{k\to \infty} \frac{1}{k} \nu_x(\til{g}^k)$$ on the basepoint $x$. Take two basepoints $x$ and $x'$ and let $\{x_s\}_{s=0}^1$, $x_0 = x$, $x_1 = x'$ be a path in $\X$ connecting them. Note that it is enough for us to show that $\nu_{x}$ and $\nu_{x'}$ differ by a bounded function $\til{\G} \to \R$. Let us compare $\nu_x(\arr{g})$ and $\nu_{x'}(\arr{g})$. Let $\delta:=\arr{g} \cdot x \#\; [g \cdot x, x]$, $\delta':=\arr{g} \cdot x' \#\; [g\cdot x', x']$ and let $D,D'$ be their contracting discs. Define the disk $C:[0,1]\times [0,1] \to \X$ by $C(s,t)= g_t \cdot x_s$. Moreover define $S^0$ be the contracting disk of $\{x_s\} \#\; [x',x]$. Then by preservation of $\K$ by the action $S^1= g \cdot S^0$ will be the contracting disk of $\{g \cdot x_s\} \# [g\cdot x', g \cdot x]$. Note that $g \cdot x_s \equiv C(s,1)$. At last, define an adapted contracting disk of $[x,x']\cup[x,g \cdot x']\cup[g \cdot x, x']\cup[g \cdot g, g \cdot x']$ as the union $Q = \Delta_0 \cup \Delta_1$ for the two geodesic triangles $\Delta_0 , \Delta_1$ on $\{x,x',g \cdot x\}$ and on $\{g \cdot x,x', g \cdot x'\}$. Note then that $\Sigma = \overline{D_0}\cup D_1 \cup \overline{S^0} \cup S^1 \cup C \cup Q$ where the gluings go along the overlapping paths, is a sphere. Therefore \begin{equation*}  0 = \int_\Sigma \Om = \int_{\overline{D_0}\cup D_1 \cup \overline{S^0} \cup S^1 \cup C \cup Q} \Om = \int_{D_1} \Om - \int_{D_0} \Om - \int_{S^0}\Om + \int_{g \cdot S^0} \Om + \int_Q \Om + \int_C \Om = \end{equation*} \begin{equation}\label{symplectic area part of independence upon endpoint} = \int_{D_1} \Om - \int_{D_0} \Om  + \int_Q \Om + \int_C \Om = \nu_x(\til{\phi}) -\nu_{x'}(\til{\phi}) + \intoi \mu(X_t)(g_t \cdot x) - \intoi \mu(X_t)(g_t \cdot x') + \int_C \Om  + \int_Q \Om\end{equation} since the action of $\til{\G}$ on $\X$ preserves $\Om$. Wherefrom \begin{equation}\label{indep on basepoint with Z}|\nu_x(\til{g}) -\nu_{x'}(\til{g})| \leq |Z| + 2C_\X,\end{equation} for $Z=\intoi \mu(X_t)(g_t \cdot x) - \intoi \mu(X_t)(g_t \cdot x') + \int_C \Om $. We now show that $Z$ equals zero, finishing the proof. And indeed letting $\paph{X_t}$ be the path in $Lie(\G)$ corresponding to $\arr{g}$ and $\Xi_t = \overline{X}_t$, we have

$$\int_C \Omega = \int_0^1 \int_0^1 \Om(\del_s C(s,t),\del_t C(s,t)) ds dt = \int_0^1 \int_0^1 \Om(\del_s C,\Xi_t(C(s,t))) ds dt =$$$$= \int_0^1 \int_0^1 d_{C(s,t)}\mu(X_t)(\del_s C) ds dt = \int_0^1 \int_0^1 \del_s \mu(X_t)(C(s,t)) ds dt =$$$$= \int_0^1 \mu(X_t)(C(1,t)) - \mu(X_t)(C(0,t)) dt = \int_0^1\mu(X_t)(g_t \cdot x')\,dt - \int_0^1\mu(X_t)(g_t \cdot x)dt.$$
\end{proof}

We also prove Proposition \ref{general symp independence} on the transformation of $\nu_x$ under conjugation with respect to a suitable normal extension.

\begin{proof}
Consider a path $\arr{g}=\paph{g_t}$ representing $\til{g} \in \til{\G}$. Then for an element $h \in \cH$ the path $h\arr{g}h^{-1} = \paph{h g_t h^{-1}}$  will represent $h\til{g}h^{-1}$. By definition \[\nu_x(h\arr{g}h^{-1}) = \int_{D_{h\arr{g}h^{-1}}} \Om - \intoi \mu(Ad_h X_t)(h g_t h^{-1}\cdot x) dt =\] for a disk $D_{h\arr{g}h^{-1}}$ with boundary $\delta_x=\arr{g}\cdot x \#\; [g\cdot x, x]$, and noting that by preservation of $\K$ we have the relation  $h\cdot \delta_{h^{-1}\cdot x} = \delta_x$ for $\delta_{h^{-1}\cdot x} = \arr{g}\cdot (h^{-1}\cdot x) \#\; [g\cdot h^{-1}\cdot x, h^{-1}\cdot x]$ so that the disk $D$ satisfying \[h\cdot D= D_{h\arr{g}h^{-1}}\] has boundary $\delta_{h^{-1}\cdot x}$, so that by preservation of $\Om$ and by equivariance of $\mu$ with respect to $\cH$ we have
\[ = \int_{D} \Om - \intoi \mu(X_t)(g_t h^{-1}\cdot x) dt = \nu_{h^{-1}\cdot x}(\arr{g}),\] which proves the proposition. Note that for every $\til{h} \in \til{\cH}$ with endpoint $h$ we have $\til{h} \til{g} \til{h}^{-1} = h \til{g} h^{-1}$ since the paths $\paph{h_t g_t h_t^{-1}}$ and $h\arr{g}h^{-1}$ are homotopic with fixed endpoints.
\end{proof}

\subsection{Finite dimensional examples and Guichardet-Wigner quasimorphisms}\label{finite dimensional and Guichardet-Wigner}

In this section we define the Guichardet-Wigner quasimorphisms and prove Proposition \ref{finite dimensional nu and varrho GW} on reconstructing these through moment maps.

We remark that as we have assumed that $G$ has finite center, there are no homogenous quasimorphisms on $G$ (cf. \cite{ShternAutomaticContinuity,BurgerMonodLattices} and \cite{BargeGhysEulerMaslov} for the group $Sp(2n,\R)$). Moreover it is known that the all homogenous quasimorphisms on $\til{G}$ are proportional to $\varrho_G$ (cf. \cite{ShternAutomaticContinuity,BurgerMonodLattices,BargeGhysEulerMaslov}). From these two remarks it follows that it is enough to show the equality of $\nu_G$ and $\varrho_{G}$ on $\pi_1(G) \cong \pi_1(K)$. In fact, $\varrho_{G}$ is defined as the unique homogenous quasimorphism $\til{G} \to \R$ such that its pullback $\varrho_{G}|_{\til{K}}: \til{K} \to \R$ to $\til{K}$ by the natural map $\til{K} \to \til{G}$ coincides with the lift $\til{v}:\til{K} \to \R$ of the canonical (up to powers) character $v: K \to S^1$, constructed in either one of several ways. The first way is as follows.
The Lie algebra $\mathfrak{k}$ of $K$ satisfies $\fk= \fz + [\fk,\fk]$ where $\fz$ is the center of $\fk$ (Corollaries 4.25 and 1.56 in \cite{Knapp}). In the case when $G$ is a simple Hermitian symmetric Lie group, $\fz$ is one-dimensional by \cite{Knapp} p.513. Hence the center $Z$ of $K$ is one dimensional. Take the identity component $Z_0 \cong S^1$  of $Z$. Then by Theorem 4.29 in \cite{Knapp} $K=(Z_0)K_{ss}$, for $K_{ss}$ the analytic subgroup with Lie algebra $[\fk,\fk]$. The group $K_{ss}$ has a finite center, therefore by taking quotients by $K_{ss}$ we get a homomorphism $v: K \to Q \cong S^1$ from $K$ to the quotient $Q \cong S^1$ of $Z_0 \cong S^1$ by a finite subgroup.

\begin{example}\label{first way for G=Sp(2n,R)}
For $G=Sp(2n,\R)$ we have $K\cong U(n)$ and $K_{ss} \cong SU(n)$. Therefore the first construction gives the homomorphism $v: K \to U(n)/SU(n) \cong S^1$ is simply $v(k) = det_\C(k)$.
\end{example}

The second way to construct $v: K \to S^1$ is by use of the action of $G$ on the Hermitian symmetric space $X = G/K$ - it is shown in \cite{GuichardetWigner} that $v$ equals the determinant of the linearization of the natural action of $K \subset G$ at the fixed point $x=[Id] \in X = G/K$. Note that the two constructions of $v$ agree up to the power $-2\dim_{\C}(X)/\#(Z_0 \cap K_{ss})$ since the determinant of a scalar matrix equals the scalar raised to the power of the dimension of the space (cf. \cite{Helgason} - proof of Theorem 6.1 and \cite{GuichardetWigner} - proof of Th$\acute{\text{e}}$or$\grave{\text{e}}$me 2).

\begin{example}\label{second way for G=Sp(2n,R)}
For $G=Sp(2n,\R)$ we have $K_{ss} \cong SU(n)$, $Z_0 \cong D$ - the subgroup of diagonal matrices in $U(n)$ and $\#(Z_0 \cap K_{ss}) = n$. As in this case $\dim_\C(X) = n(n+1)/2$, the second construction gives the homomorphism $v(k) = det_\C^{-(n+1)}(k)$.
\end{example}

We use the second way to define $\rho_G$ now. Proposition \ref{finite dimensional nu and varrho GW} is then demonstrated as follows.

\begin{proof}

Consider the point $x=[Id] \in X = G/K$. It is a fixed point under the natural action of $K \subset G$. By the construction of $\nu_x$ and of the equivariant moment map $\mu: Lie(G) \times X \to \R$ for a path $\arr{k}=\paph{k_t}$ in $K$ with $k_0 = Id$ representing $\til{k} \in \til{K}$, we have \[\nu_x(\til{k}) = -\intoi \mu(\eta_t)(x) dt = -\intoi (\alpha_X)_y(\widehat{\eta}_y) dt = -\frac{1}{i} \intoi \frac{d}{dt'}|_{t'=t} det_{\C}((k_{t'})_{*_x} ) det_{\C}((k_{t})_{*_x} )^{-1}dt = .\] \[= -\text{varangle}(\paph{det_{\C}((k_{t})_{*_x})}) = -\til{v}(\til{k}).\]
Hence $\nu_x$ equals $\til{v}$ on $\til{K}$, and consequently $\nu_x$ equals $\varrho_G$ on $\pi_1(G) \cong \pi_1(K)$. Therefore the homogeneization $\nu_G$ of $\nu_x$ equals $-\varrho_G$ on $\pi_1(G)$ and this confers the equality $\nu_G = -\varrho_G$ on the whole group $\til{G}$.

\end{proof}

\subsection{The equality of the homomorphisms $\A$ and $I_{c_1}$ on $\pi_1(Ham(M,\om))$}
Now we prove Theorem~\ref{Theorem Equality of A_mu and I_c1}.

\begin{proof}
First we note the following equality due to Fujiki \cite{Fujiki}.
Given the bundle $\mathcal{Z}$ over $\J$ which has $\ZZ_J:=(M,J)$ for the fiber over $J$, denote by $T_{\ZZ,\J}$ the vertical bundle and take $c_1(K)$ to be the Chern form of the vertical canonical bundle $K$ relative to the Hermitian metric given by $h(J) = g(J) - i \om$ in the fiber over $J \in \J$, then \begin{equation}\label{FujikiFormula}\Om = \int_{fiber} c_1(K) \; p^*\om^n,\end{equation} for $p:\mathcal{Z} \to M$ the smooth projection map.

%

The Hamiltonian fiber bundle over $S^2$ corresponding to a loop $\gamma = \paph{\phi_t}$ in $\G$ based at the identity can be described (cf. \cite{PolterovichLoops}) as $P_\gamma = M \times D_{-} \cup_{\Phi} M \times D_{+}$, where $D_{-}$ and $D_{+}$ are two copies of the disk $\D$ and the gluing map $\Phi:\del(M\times D_{-}) \cong M \times S^1 \to M \times S^1 \cong \del(M \times D_{+})$ is given by $\Phi:(x,t) \mapsto (\phi_t x, t)$.

Note that given a Hamiltonian loop $\gamma = \{\phi_t\}_{0 \leq t \leq 1}$ the bundle $P=P_\gamma$ with a vertical compatible complex structure is obtained by a map $D:\D \to \J$ representing a relative homotopy class in $\pi_2(\J,\G_{J_0})$ corresponding to the loop $\gamma^{-1} = \{\psi_t = (\phi_t)^{-1}\}_{0 \leq t \leq 1}$ - that is $\del D : S^1 \to \J$ is given by $\paph{(\psi_t)_*J_0}$. Note that $P|_{D_-}$ with its fiberwise complex structure is equal to $D^*\ZZ$. We denote by $H_t$ the zero-mean-normalized Hamiltonian for $\gamma$ and by $G_t$ the zero-mean-normalized Hamiltonian for $\gamma^{-1}$. The two are related by the formula $G_t(x) = - H_t(\phi_t x)$.

Moreover $$-I_{c_1}(\gamma) = \int_{\D} \int_{fiber} D^*c_1(K) u^n.$$ Since the coupling class $u$ is represented by the form $\Upsilon :=$ \{$\omega$ on $M \times D_{+}$; $\omega + d (\psi(r) H_t(\phi_t x) dt))$ on $M \times D_{-}$ \}, we have $$-I_{c_1} = \int_{\D} \int_{fiber} D^*c_1(K) (\omega^n + n d (\psi(r) H_t(\phi_t x) \, dt \, \omega^{n-1}))$$ $$= \int_{\D} \int_{fiber} D^*c_1(K) \omega^n + n \int_{D} \int_{fiber} D^*c_1(K) d (\psi(r) H_t(\phi_t x) \, dt \, \omega^{n-1}).$$ By the result of Fujiki the first summand equals $\int_\D D^*\Om$. It is therefore enough to show that the second summand equals $\int_0^1 dt \int_M S(\psi_t \cdot J_0)G_t(x)\om^n$. The second summand satisfies $$ n \int_{\D} \int_{fiber} D^*c_1(K) d (\psi(r) H_t(\phi_t x) \, dt \, \omega^{n-1}) = n \int_{M \times \D} d ( D^*c_1(K)\psi(r)H_t(\phi_t x) dt \om^{n-1}) = $$

$$ = n \int_{M \times S^1} H_t(\phi_t x) D^*c_1 \om^{n-1} dt =$$ and by Equation \ref{definition of S(J)} we have


$$= \intoi \int_M S(\psi_t\cdot J_0) H_t(\phi_t x) \om^n(x)  dt =  -\intoi \int_M S(\psi_t\cdot J_0) G_t (x) \om^n(x)  dt.$$

Consequently we have $I_{c_1}(\gamma)= -\A_\mu(\gamma^{-1}) = A(\gamma)$.

\end{proof}

\subsection{The finite-dimensional case $G=Sp(2n,\R)$ and the Maslov quasimorphism}\label{FinDimExamplesSp}

In this section we would like to write out the finite-dimensional example more explicitly in the case $G=Sp(2n,\R)$ - for later use in particular. When $G=Sp(2n,\R)$ the maximal compact subgroup is $K \cong U(n)$ and the space $X=G/K$ has several guises. First it can be considered as the Siegel upper half-space \cite{Siegel} $\cS_n = \{X+iY|\;X,Y \in \text{Mat}(n,\R),\,X=X^t,Y=Y^t,\, Y>0\} \subset \text{Mat}(n,\C)$. Here there is a natural Kahler form $\sigma_{Siegel} = \trace(Y^{-1}dX \wedge Y^{-1}dY)$ where the complex structure comes from the one on $\text{Mat}(n,\C)$. This form is Kahler-Einstein with cosmological constant $\lambda = -\frac{n+1}{2}$ \cite{Siegel} - that is \begin{equation}\label{sigma_Siegel cosmological constant} Ric(\sigma_{Siegel}) = -\frac{n+1}{2} \sigma_{Siegel}.\end{equation} From which, since proportional metrics have equal Ricci forms, we have \begin{equation}\label{sigma_Siegel and sigma_Bergman}\sigma_{Siegel} = \frac{2}{n+1}\sigma_{Bergman},\end{equation} for the Bergman Kahler form $\sigma_{Bergman}$ on $X$.

Second, the space $X=G/K$ can be considered as the space $\rJ_c$ of $\om_{std}$-compatible complex structures on the symplectic vector space $(\R^{2n},\om_{std})$. In this model, a natural symplectic form $\sigma_{\trace}$ is given by $(\sigma_{\trace})_J(A,B) = \frac{1}{4}\trace(JAB)$ for $J \in \rJ_c$ and $A,B \in T_J (\rJ_c)$. A short computation based on the fact that all $G$-invariant $2$-forms on $X$ are proportional gives \begin{equation}\label{sigma_trace and sigma_Siegel}\sigma_{\trace} = \frac{1}{2} \sigma_{Siegel},\end{equation} under the natural isomorphisms $\rJ_c \cong X$, $\cS_n \cong X$.

By Examples \ref{first way for G=Sp(2n,R)} and \ref{second way for G=Sp(2n,R)} the Maslov quasimorphism $\tau_{Lin}: \til{G} \to \R$ restricting on $\til{K}$ to $\til{v}$ for $v= det_\C^2$ on $K \cong U(n)$ can be written as \begin{equation}\label{tau_Lin via nu_G Bergman}\tau_{Lin} = \frac{2}{n+1} \nu_{G,Bergman}\end{equation} in terms of $\nu_G$ for $\sigma_X =\sigma_{Bergman}$. Therefore by Equation \ref{sigma_Siegel and sigma_Bergman} \begin{equation}\label{tau_Lin via nu_G Siegel}\tau_{Lin} = \nu_{G,Siegel}\end{equation} for $\sigma_X =\sigma_{Siegel}$, and by Equation \ref{sigma_trace and sigma_Siegel} \begin{equation}\label{tau_Lin via nu_G trace}\frac{1}{2}\tau_{Lin} = \nu_{G,\trace}\end{equation} for $\sigma_X =\sigma_{\trace}$.

Note that by \cite{Siegel} $\sigma_{Siegel}$ and consequently $\sigma_{\trace}$ has non-positive sectional curvature. Moreover $\sigma_{\trace}$ is Kahler-Einstein with cosmological constant $-(n+1)$.

Now consider $X \cong \rJ_c$ with $\sigma_X = \sigma_{\trace}$. By Equation \ref{tau_Lin via nu_G trace}, Lemma \ref{nu of x and nu of x^{-1}} and by the definition of $\nu_x$ we have $-\frac{1}{2}\tau_{Lin} = \nu_G \simeq \nu_x(\arr{g}) =\int_{D_{\arr{g}}} \sigma_X - \intoi \mu(\xi_t)(g_t\cdot x) dt$. Hence \begin{equation}\label{int_D sigma via tau_Lin and mu}\int_{D_{\arr{g}}} \sigma_X \simeq \frac{1}{2}\tau_{Lin}(\arr{g}) + \intoi \mu(\xi_t)(g_t\cdot x) dt.\end{equation}

 For later calculations we will want the moment map summand in this formula more explicit. We write a formula for $\mu$ using the fact that it is an equivariant moment map for the action of the semisimple Lie group $Sp(2n,\R)$ (cf. \cite{Knapp}) on $\cS_n$. Note that an equivariant moment map for the action of a semisimple Lie group on a symplectic manifold is unique \cite{IntroSymp}. Hence it is enough to show the following.

\begin{lma}\label{equivariant moment map for the linear symplectic group}
Consider the action of $Sp(2n,\R)$ on $\cS_n \cong \rJ_c$ with the invariant Kahler form $\sigma_{trace}$. Then it is Hamiltonian with the equivariant moment map $\mu_{\cS_n}: \cS_n \times \mathfrak{sp}(2n,\R) \to \R$ given by $\mu_{\cS_n}(J)(\Xi) = - \frac{1}{2} trace(\Xi J)$.
\end{lma}

\begin{pf}
The symplectic form $\sigma$ on $\cS_n$ can be described as $\sigma(A,B) = \frac{1}{4} \trace(JAB)$ using the isomorphism $\cS_n \cong \rJ_c$ - the space of complex structures on $\R^{2n}$ compatible with the standard symplectic form. Let us first compute the vector field $\overline{\Xi}$ generated by the infinitesimal action of $\Xi$. At a point $J \in \cS_n$, denoting $\Phi_t = exp(t\Xi) \in Sp(2n,\R)$ we have $\overline{\Xi}_J = \frac{d}{dt}|_{t=0} \Phi_t \cdot J = \frac{d}{dt}|_{t=0} \Phi_t J \Phi_t^{-1} = \Xi J - J \Xi = -[J,\Xi]$. Then for $B \in T_J \cS_n$ we compute $$d_J (\trace(\Xi J))(B) = \trace(\Xi B).$$ Finally, for $B \in T_J \cS_n$ we have $$(i_{\overbar{\Xi}}\sigma)_J(B) = \sigma_J(\overline{\Xi}_J,B) = -\sigma_J([J,\Xi],B) = -\frac{1}{4}\trace(J[J,\Xi]B) =$$ $$= -\frac{1}{4}\trace(-J\Xi JB + J^2 \Xi B) = \frac{1}{4}\trace(\Xi JBJ + \Xi B) = \frac{1}{2} \trace( \Xi B).$$ The last expression equals $- d_J (-\frac{1}{2}\trace(\Xi J))(B)$ as we have computed, and we're done.
\end{pf}

\subsection{The local type of the quasimorphism on $\til{Ham}(M,\om)$}

We shall now describe the local behaviour of the quasimorphism $\fS$ - we compute its restriction to subgroups $\G_B \subset \G$ of diffeomorphisms supported in embedded balls $B$ in $M$, proving Theorem \ref{Theorem local type}.

\begin{df}\label{small ball}(Embedded balls) We denote by $\cU$ the set of embedded balls in $M$.
\end{df}

Given a symplectic manifold $(M,\om)$ with an almost complex structure $J_0 \in \J$ with Hermitian scalar curvature $S(J_0)$ of mean $c = n{\int_M c_1(TM,\om)\om^{n-1}}/{\int_M \om^n}$, and $B \in \cU$ an embedded ball in $M$, we will show that the restriction $\nu_B=\fS|_{\G_B}$ of the quasimorphism $\fS$ to $\G_B = Ham_c(B,\om|_B)$ satisfies $$\nu_B = \frac{1}{2}\tau_B - c \textit{Cal},$$ where $\tau_B$ is the Barge-Ghys Maslov quasimorphism on $\G_B = Ham_c(B^{2n},\om_{std})$ and $\textit{Cal}$ is the Calabi homomorphism.

Since the quasimorphism $\fS$ is homogenous and its distance from $\fS_{J_0}$ is bounded we can make calculations with $\fS_{J_0}$ allowing for an error term that vanishes under homogenization. The proof consists of writing $\nu_B$ (using Section \ref{FinDimExamplesSp}) as the sum of $\frac{1}{2}\tau$  and a remainder term. Then we use some differential geometry to show that the remainder term equals a multiple of the Calabi homomorphism. For the differential geometry part we would like to use the canonical connection on the Hermitian manifold $(M,\om,J_0)$ that is defined by the following of its properties. It preserves $\om$ and $J_0$ and its torsion has vanishing $(1,1)$-component: \begin{equation}\label{properties of Chern connection generalized}\nabla J_0 = 0, \nabla \om = 0, T^{(1,1)}_\nabla \equiv 0.\end{equation}
This connection has an equivalent definition in terms of $\bar{\partial}$-operators on complex vector bundles, which is the one used in \cite{DonaldsonGaugeTheoryComplexGeometry}. It is sometimes called "the Chern connection", and sometimes - "the second canonical connection of Ehresmann-Liebermann" (cf. \cite{Gauduchon} Section 2, \cite{KobayashiConnectionsAlmsotComplex} and \cite{TosattiWeinkoveYauTamingSymplecticForms} Section 2).

Consider $B$ as a smooth embedding $B:B^{2n} \to M$ from the standard ball \[B^{2n} = \{(z_1,...,z_n)|\Sum_{j=1}^n |z_j|^2 < 1\} \subset \C^n\] to $M$. For purposes of trivialization and estimates choose for each two points $x,y \in B$ a path $\gamma_{x,y}$ starting at $x$ and ending at $y$, that depends continuously on $(x,y) \in B\times B$ where $\gamma_{x,x}$ is the constant path at $x$ for all $x \in B$. This can be achieved for example by taking linear segments in $B^{2n}$. Then we have the following lemma.

\begin{lma}\label{bounded and continuous trivializations}
Let $B_{-}\Subset B$ be any closed ball compactly contained in $B$. Then the following two statements hold by continuity and compactness of $B_{-}\times B_{-}$.
\begin{enumerate}
\item For every one-form $\lambda \in \Om^1(B)$ the function $B_{-}\times B_{-} \to \R$ given by $(x,y) \mapsto \int_{\gamma_{x,y}}\lambda$ is bounded by a constant depending only on $B,B_{-},\lambda$.
\item Given any connection $\nabla'$ preserving $\om$ and any fixed symplectic trivialization $TB \cong V \times B$ for a symplectic vector space $(V,\om_V)$, the map $B_{-} \times B_{-} \to Sp(V,\om_V)$ obtained by means of the trivialization by the parallel transport $\Gamma_{\gamma_{x,y}}:T_x B \to T_y B$ with respect to $\nabla'$ has a compact image in $Sp(V,\om_V)$.
\end{enumerate}
\end{lma}

Take a path $\paph{\phi_t} \subset \G_B$ with $\phi_0 = Id$. We shall now unwind the definition of $\nu_B(\paph{\phi_t})$. Over each $x \in B$ we have the fiber $\cS_x$ of the bundle $\cS \to B$. In $S_x$ we have the path $(\phi_t\cdot J_0)_x$. Now we shall define a path $\Phi(x)_t$ in $Sp(T_xM,\om_x)$ associated to ${(\phi_t)_*}_x$ such that under the action of $Sp(T_xM,\om_x)$ on $S_x$, we have  $\Phi(x)_t\cdot (J_0)_x = (\phi_t\cdot J_0)_x = ({{\phi_t}_*}_{\phi_t^{-1}x}) (J_0)_{\phi_t^{-1}x} ({{\phi_t}_*}_{\phi_t^{-1}x})^{-1}$.

Indeed consider for each $t \in [0,1]$ the path $\gamma_{x,\phi_t^{-1} x}$. The parallel transport along this path preserves $J_0$ and maps $\Gamma_{\gamma_{x, \phi_t^{-1} x}}: T_x M \to T_{\phi_t^{-1} x} M$. Then $\Phi(x)_t=   ({{\phi_t}_*}_{\phi_t^{-1}x}) \circ\Gamma_{\gamma_{x, \phi_t^{-1} x}} : T_x M \to T_x M$ is the required map. Indeed $\Phi(x)_t \cdot (J_0)_x = \Phi(x)_t (J_0)_x \Phi(x)_t^{-1} =({{\phi_t}_*}_{\phi_t^{-1}x}) \Gamma_{\gamma_{x, \phi_t^{-1} x}} (J_0)_x  (\Gamma_{\gamma_{x, \phi_t^{-1} x}})^{-1} ({{\phi_t}_*}_{\phi_t^{-1}x})^{-1}= ({{\phi_t}_*}_{\phi_t^{-1}x})  (J_0)_{\phi_t^{-1} x}  ({{\phi_t}_*}_{\phi_t^{-1}x})^{-1} = (\phi_t \cdot J_0)_x$, by preservation of $J_0$. Henceforth we omit the subscript $z$ in $(J_0)_z$ whenever this is determined by the context.

Then for all $x\in B$ we have the loop $\delta(x) = \paph{\Phi(x)_t \cdot J_0} \# \overline{[J_0, \Phi(x)_1 \cdot J_0]}$. We then for all $x \in B$ choose a disk $D(x)$ that bounds $\delta(x)$ - in fact one can construct $D(x)$ as the geodesic join of $\paph{\Phi(x)_t \cdot J_0}$ with $J_0$ - that is $D(x) = \bigcup_t \overline{[J_0, \Phi(x)_t \cdot J_0]}$ properly parametrized. Denote $\gamma_t(x)= \overline{[J_0, \Phi(x)_t \cdot J_0]}$. Denote by $\beta_t(x)$ the path $\{\Phi(x)_{t'} \cdot J_0\}_{t'=0}^{t'=t}$.

Recall from Section \ref{preliminaries on quasimorphisms} that $a \simeq b$ denotes the equality of the functions $a,b$ up to a function that is bounded by a constant that does not depend on their arguments. Compute \begin{equation}\label{write out nu on a ball}\nu_B(\paph{\phi_t}) \simeq \int_D \Om - \intoi \mu(X_t)(\phi_t \cdot J_0) = \int_B (\int_{D(x)} \sigma_x) \om^n(x) - \intoi \int_M S(\phi_t \cdot J_0) H_t(x) \om^n(x).\end{equation} Now note that by Equation \ref{int_D sigma via tau_Lin and mu} and the definition of the moment map for the action of $G=Sp(2n)$ on $X=G/K$ \begin{equation}\label{local at x}\int_{D(x)}\sigma_x \simeq \frac{1}{2}\tau_{Lin}(\paph{\Phi(x)_t}) + \int_0^1 h(x)_t(\Phi(x)_t \cdot (J_0)_x) dt, \end{equation} where the function $h(x)_t(\cdot) = \mu_{S_x}(\Xi(x)_t)(\cdot)$, for $\Xi(x)_t = \frac{d}{d\tau}|_{\tau = t}\Phi(x)_\tau \circ \Phi(x)_t^{-1}$, is the contact Hamiltonian for the canonical lifting of $\Phi(x)_t$ to the principal $S^1$-bundle of unit vectors in $\Lambda^N T\cS_x$ simply by use of the differential (cf. Equation \ref{equation - moment map for finite dimensional groups}). As a side remark it may be said, following \cite{DonaldsonMomentMapsDiffGeom}, that this finite-dimensional moment map is the main reason for the existence of the corresponding infinite-dimensional moment map.

 Consequently, integrating over $B$ with respect to the form $\om^n$, we have \begin{equation} \nu_B(\paph{\phi_t}) \simeq \frac{1}{2} \cdot \int_B \tau_{Lin} (\paph{{\Phi(x)_t}})\om^n(x) + \int_B \int_0^1 h(x)_t(\Phi(x)_t \cdot J_0) \, dt \, \om^n(x) - \intoi \int_M S(\phi_t \cdot J_0) H_t(x) \om^n(x).\end{equation} By the definition of the Barge-Ghys Maslov quasimorphism on $\G$ and Lemma \ref{bounded and continuous trivializations}, the first term homogenizes to $\frac{1}{2}\tau_B$. Our goal is hence to compute the sum of the second and the third terms.

 By Lemma \ref{equivariant moment map for the linear symplectic group} we rewrite the second term in Equation \ref{local at x} as \begin{equation}\label{rewrite local at x with trace}\int_0^1 h(x)_t(\Phi(x)_t \cdot J_0) = -\int_0^1 \frac{1}{2} trace(\Xi(x)_t(\Phi(x)_t \cdot J_0)) dt ,\end{equation} for $\Xi(x)_t = \frac{d}{d\tau}|_{\tau = t}\Phi(x)_\tau \circ \Phi(x)_t^{-1}$.

Now note that instead of using the parallel transport along $\gamma_{x,\phi_t^{-1} x}$ to define $\Phi(x)_t :T_x B \to T_x B$ we could use the one along $p_{x,t} = \{\phi_{t'}^{-1} x\}_{t'=0}^t$ to define the map $\Psi(x)_t = ({{\phi_t}_*}_{\phi_t^{-1}x}) \circ\Gamma_{p_{x,t}}: T_x B \to T_x B$. Then we have \begin{equation}\label{Phi via Psi and U}\Phi(x)_t = \Psi(x)_t U(x,t),\end{equation} for the unitary map $U(x,t) = \Gamma_{p_{x,t}}^{-1} \circ \Gamma_{\gamma_{x,\phi_t^{-1} x}} : T_x B \to T_x B$. Form $\Upsilon(x)_t = \frac{d}{d\tau}|_{\tau = t}\Psi(x)_\tau \circ \Psi(x)_t^{-1}$ and $\Theta(x,t) = \frac{d}{d\tau}|_{\tau = t}U(x,\tau) \circ U(x,t)^{-1}.$ Then by Equation \ref{Phi via Psi and U} we have \begin{equation}\label{Xi via Theta and Upsilon}\Xi(x)_t = \Psi(x)_t\Theta(x,t)\Psi(x)_t^{-1} + \Upsilon(x)_t,\end{equation} and \begin{equation}\label{Phi action via Psi action and U} \Phi(x)_t \cdot J_0 = \Phi(x)_t J_0 \Phi(x)_t^{-1} =   \Psi(x)_t U(x,t) J_0 U(x,t)^{-1} {\Psi(x)_t}^{-1} = \Psi(x)_t  J_0 {\Psi(x)_t}^{-1}=\Psi(x)_t \cdot J_0,\end{equation} because $U(x,t)$ is $J_0$-linear.

Therefore, by Equation (\ref{rewrite local at x with trace}) and noting that $$\trace(\Psi(x)_t\Theta(x,t)\Psi(x)_t^{-1}(\Psi(x)_t \cdot J_0)) = \trace(\Theta(x,t)J_0)$$ we have \begin{equation}\label{rewrite local at x with Theta and Upsilon}\int_0^1 h(x)_t(\Phi(x)_t \cdot J_0) = -\int_0^1 \frac{1}{2} \trace(\Theta(x,t)\cdot J_0) dt -  \int_0^1 \frac{1}{2} \trace(\Upsilon(x)_t(\Psi(x)_t \cdot J_0)) dt.\end{equation}

Note additionally, that $$\frac{1}{2}\trace(\Theta(x,t)\cdot J_0) = -\frac{1}{i}\trace_\C(\Theta(x,t)),$$ considering $\Theta(x,t)$ as a skew-Hermitian operator on the complex Hermitian space $(T_x B, J_0, \om_x)$. Moreover, $$\trace_\C(\Theta(x,t)) = \Theta^n(x,t),$$ where $\Theta^n(x,t) = \frac{d}{d\tau}|_{\tau = t}U^n(x,\tau) \circ {U^n(x,t)}^{-1},$ for $U^n(x,t) = (\Gamma^n_{p_{x,t}})^{-1} \circ {\Gamma^n_{\gamma_{x,\phi_t^{-1} x}}} : {\Lambda^n_\C} T_x B \to {\Lambda^n_\C} T_x B,$ for the naturally induced parallel translations on the Hermitian complex line bundle ${\Lambda^n_\C} TB$, endowing $TB$ with the Hermitian structure $(J_0,\om)$ and the connection $\nabla$.

Therefore \begin{equation}\label{local at x with rho}\int_0^1 h(x)_t(\Phi(x)_t \cdot J_0) = \int_{D_B(x)} \rho - \frac{1}{2} \int_0^1\trace(\Upsilon(x)_t (\Psi(x)_t \cdot J_0)),\end{equation}
where $i \rho$ is the curvature two-form of the connection $\nabla^n$ on $\Lambda^n_\C(TM,J_0)$ naturally induced from $\nabla$ on $(TM,J_0)$ and $D_B(x)$ is the disk spanned by $\bigcup_{t=0}^1 \gamma_{x,\phi_t^{-1} x}.$ Note that $\del D_B(x) = p_{x,1} \# \overline{\gamma_{x,\phi_1^{-1} x}}$.
Now $\rho|_B \in \Om^2_{closed}(B,\R)$ has by the Poincar$\acute{\text{e}}$ lemma a primitive $\alpha \in \Om^1(B,\R).$ Hence by Stokes' formula we have \begin{equation}\label{from rho to alpha on the boundary} \int_{D_B(x)} \rho = \int_{p_{x,1}}\alpha - \int_{\gamma_{x,\phi_1^{-1} x}} \alpha.\end{equation}

Choosing $B_{-} \Subset B$ such that $supp(\phi_t) \subset B_{-}$ for all $t \in [0,1]$, we have $\gamma_{y,\phi_1^{-1} y} \equiv y$ for all $y \in B\setminus B_{-}$, hence by Lemma \ref{bounded and continuous trivializations} we have the following uniform estimate for the second term $|\int_{\gamma_{x,\phi_1^{-1} x}} \alpha| \leq C(B_{-},B),$ for a constant $C(B_{-},B,\alpha)$ depending only on $\alpha$, $B_{-} \supset \bigcup_{t=0}^1 supp(\phi_t)$ and on $B$.

Now denote $\psi_t = \phi_t^{-1}$. Denote $Y_t$ the Hamiltonian vector generating $\psi_t$. Recall that $p_{x,1} = \{\phi_t^{-1} x\}_{t=0}^1 = \{\psi_t x\}_{t=0}^1.$ Hence the first term in Equation \ref{from rho to alpha on the boundary} satisfies \[\int_{p_{x,1}}\alpha = \int_0^1 ((\psi_t)^* i_{Y_t} \alpha)_x dt.\]

Hence integrating Equation \ref{local at x with rho} over $B$ we express  $$\int_B \int_0^1 h(x)_t(\Phi(x)_t \cdot J_0) \, dt \, \om^n(x)=$$ as \begin{equation}\label{writing out integral over the ball of the contact Hamiltonian}=  \int_0^1 \int_B (\psi_t)^* i_{Y_t} \alpha \, \om^n dt - \frac{1}{2} \int_B \int_0^1\trace(\Upsilon(x)_t (\Psi(x)_t \cdot J_0)) dt \om^n(x) + Bdd(\{\phi_t\}_{t=0}^1)
\end{equation}
for a function $Bdd(\{\phi_t\}_{t=0}^1)$ that satisfies $$|Bdd(\{\phi_t\}_{t=0}^1)| \leq C_1(B_{-},B,\alpha)$$ for a constant $C_1(B_{-},B)$ depending only on $\alpha$, $B_{-} \supset \bigcup_{t=0}^1 supp(\phi_t)$ and on $B$.

We shall now show that the first term in Equation \ref{writing out integral over the ball of the contact Hamiltonian} corrected by the moment map term $\intoi \mu(X_t)(\phi_t \cdot J_0) dt$ in the definition of the quasimorhism is proportional to the Calabi homomorphism. After that we will show that the second term vanishes.
Let  $G_t$ (for each $t \in [0,1]$) be the function that vanishes near $\del B$ and satisfies $i_{Y_t}\om = - d G_t.$ Then $i_{Y_t}\alpha \; \om^n = n \alpha \; i_{Y_t}\om \; \om^{n-1} = - n \alpha \; d G_t \; \om^{n-1}.$ Hence $$\int_B (\psi_t)^* i_{Y_t} \alpha \, \om^n = \int_B i_{Y_t} \alpha \, \om^n  = n \int_B d G_t \; \alpha \; \om^{n-1} = $$ $$ = - n \int_B G_t d\alpha \om^{n-1} = - n \int_B G_t d\alpha \om^{n-1} = - n \int_B G_t \rho \om^{n-1} = $$ and by definition of the Hermitian scalar curvature we have $$ = - \int_B G_t S(J_0) \om^{n} =$$ denoting $H^0_t$ (for each $t \in [0,1]$)  the function that vanishes near $\del B$ and satisfies $i_{X_t}\om = - d H^0_t,$ and noting that by the cocycle formula \cite{PolterovichGeometryGroupSymp} $G_t(x) = - H^0_t(\phi_t x),$ we have $$ = \int_B S(J_0) H^0_t(\phi_t x) \om^n(x).$$ Hence $$\intoi dt \int_B (\psi_t)^* i_{Y_t} \alpha \, \om^n - \intoi dt \int_M S(\phi_t \cdot J_0) H_t(x) \om^n(x) = \intoi dt \int_M S(\phi_t\cdot J_0)(H^0_t - H_t) \om^n = $$
where we extend $H^0_t$ by zero from $B$ to $M$, and noting that $H^0_t - H_t$ depends on $t$ only and equals the mean $\int_B H^0_t \om^n /\int_M \om^n$ we have \begin{equation}\label{calabi summand in nu} = - (\int_M S(\phi_t \cdot J_0) \om^n/\int_M \om^n) \intoi dt \int_B H^0_t \om^n = -c \cdot Cal_B (\{\phi_t\}_{t=0}^1).\end{equation}

Now it remains to show that $\int_0^1 \int_B \trace(\Upsilon(x)_t (\Psi(x)_t \cdot J_0)) \om^n(x) dt$ vanishes. First we would like to note that since the $(1,1)$-component of the torsion $T$ of $\nabla$ vanishes, we have \begin{equation}\label{Torsion and J_0}T(X,J_0 Y) = T(J_0 X,Y)\end{equation} for all vector fields $X,Y$ on $M$. Moreover since $\nabla$ preserves $J_0$ we have \begin{equation}\label{nabla of X and J_0} J_0 \nabla_{\bullet} X = \nabla_{\bullet} (J_0 X)\end{equation} for all vector fields $X$ on $M$, where for a vector field $Z$ on $M$, we denote by $\nabla_{\bullet} Z$ the endomorphism of $TM$ given by $Y \mapsto \nabla_{Y} Z$.

For a vector field $Z$ on $M$ define then the endomorphism $A_Z$ of $TM$ by $A_Z = L_Z - \nabla_Z$. Then by \cite{KobayashiNomizu}, Vol 1, Appendix 6, page 292, we have \begin{equation}\label{A-endomorphism via nabla and torsion}A_Z = - \nabla_{\bullet} Z - T(Z,\cdot)\end{equation} and \begin{equation}\label{divergence via A-endomorphism}- \trace A_Z = \text{div}_{\om^n}(Z),\end{equation} where $\text{div}_{\om^n}(Z) \in C^\infty_0(M,\R)$ is defined by  $$\text{div}_{\om^n}(Z)\,\om^n = L_Z \om^n.$$

Now we prove a formula relating the action of $J_0$ on $TM$ and the tensor $A_Z$. We claim that for all vector fields $X$ on $M$ we have \begin{equation}\label{J_0 and the A-endomorphism}\trace(A_X J_0) = \trace(A_{J_0 X}). \end{equation}

Indeed $$-\trace(A_X J_0) = \trace(\nabla_{\bullet}X \circ J_0 + T(X,J_0 \cdot)) = $$ by Equation \ref{A-endomorphism via nabla and torsion} $$= \trace(J_0 \nabla_{\bullet}X   + T(J_0 X, \cdot)) = $$ by Equation \ref{Torsion and J_0} $$ = \trace( \nabla_{\bullet}J_0 X  + T(J_0 X, \cdot)) = - \trace{A_{J_0 X}}, $$by Equation \ref{nabla of X and J_0}.

Let us now compute $\Upsilon(x)_t = \frac{d}{d\tau}|_{\tau = t}\Psi(x)_\tau \circ \Psi(x)_t^{-1}$ in terms of the connection and of the vector field $X_t$ generating the path of diffeomorphisms $\{\phi_t\}_{t=0}^1$. Recalling that $\Psi(x)_t = ({{\phi_t}_*}_{\phi_t^{-1}x}) \circ\Gamma_{p_{x,t}}$ we have $$\frac{d}{d\tau}|_{\tau = t}\Psi(x)_\tau =  ({{\phi_t}_*}_{\phi_t^{-1}x}) (L_{X_t} - \nabla_{X_t})_{\phi_t^{-1}(x)} \Gamma_{p_{x,t}}.$$
Consequently, \begin{equation}\label{Upsilon}\Upsilon(x)_t = ({{\phi_t}_*}_{\phi_t^{-1}x}) (A_{X_t})_{\phi_t^{-1}(x)} ({{\phi_t}_*}_{\phi_t^{-1}x})^{-1}\end{equation} for the endomorphism $A_{X_t}$ of $TM$. Then $$\trace(\Upsilon(x)_t (\Psi(x)_t \cdot J_0)) = \trace(({{\phi_t}_*}_{\phi_t^{-1}x}) (A_{X_t})_{\phi_t^{-1}(x)} ({{\phi_t}_*}_{\phi_t^{-1}x})^{-1} (({{\phi_t}_*}_{\phi_t^{-1}x})\Gamma_{p_{x,t}} (J_0)_{x} \Gamma_{p_{x,t}}^{-1} ({{\phi_t}_*}_{\phi_t^{-1}x})^{-1} )) = $$ \begin{equation}\label{trace with Upsilon} = \trace((A_{X_t})_{\phi_t^{-1}(x)} \Gamma_{p_{x,t}} (J_0)_{x} \Gamma_{p_{x,t}}^{-1}) =
\trace(A_{X_t} J_0)(\phi_t^{-1}(x)) = \trace(A_{J_0 X_t})(\phi_t^{-1}(x)),\end{equation}
by Equation \ref{J_0 and the A-endomorphism}. Hence $$\int_0^1 \int_B \trace(\Upsilon(x)_t (\Psi(x)_t \cdot J_0)) \om^n(x) dt = \int_0^1 \int_B \trace(A_{J_0 X_t})(\phi_t^{-1}(x)) \om^n(x) dt =$$\begin{equation}\label{integral of trace equals via div to zero}= \int_0^1 \int_B \trace(A_{J_0 X_t}) \om^n dt = - \int_0^1 \int_B \text{div}(J_0 X_t) \om^n dt = 0.\end{equation}

 Therefore, assembling Equations \ref{write out nu on a ball},\ref{writing out integral over the ball of the contact Hamiltonian},\ref{calabi summand in nu}, \ref{integral of trace equals via div to zero} and Definition \ref{Barge-Ghys average Maslov quasimorphism} we have \[\nu_B(\paph{\phi_t}) = \frac{1}{2} \cdot \tau_B(\paph{\phi_t}) - c \cdot Cal_B(\paph{\phi_t}) + Bdd_2(\paph{\phi_t}),\] for a function $Bdd_2(\paph{\phi_t})$ bounded by a constant $C_2(B_{-},B,\alpha)$ that depends only on $B,\alpha$ and $B_{-} \supset \bigcup_{t=0}^1 supp({\phi_t})$. Noting that $supp({\phi_t}^k) \subset supp({\phi_t})$ for every $t \in [0,1], \; k \in \Z$ and homogenizing, we finish the proof.

\subsection{The restriction to the Py quasimorphism}\label{Restriction to Py}

In this section we prove the first point of Theorem \ref{General agrees with Py and Entov} on the equality of the Py quasimorphism of Defintion \ref{Py quasimorphism} and the general quasimorphism from Corollary \ref{Ham general principle} when the symplectic manifold $(M,\om)$ is monotone - that is $c_1(TM,\om)= \kappa [\om]$ where $\kappa \neq 0$. The computation is somewhat similar to that of the local type - with the exception that there is no trivialization involved really.

 As in the computation of the local type, we use the parallel transport along $p_{x,t} = \{\phi_{t'}^{-1} x\}_{t'=0}^t$ to define the map $\Psi(x)_t = ({{\phi_t}_*}_{\phi_t^{-1}x}) \circ\Gamma_{p_{x,t}}: T_x B \to T_x B$. Then $\Upsilon(x)_t = \frac{d}{d\tau}|_{\tau = t}\Psi(x)_\tau \circ \Psi(x)_t^{-1}$ will satisfy $$\Upsilon(x)_t = ({{\phi_t}_*}_{\phi_t^{-1}x}) (A_{X_t})_{\phi_t^{-1}(x)} ({{\phi_t}_*}_{\phi_t^{-1}x})^{-1}$$ for the endomorphism $A_{X_t}$ of $TM$, for $A_{X_t} = L_{X_t} - \nabla_{X_t}$ as in Equation \ref{Upsilon}. Then $$\trace(\Upsilon(x)_t (\Psi(x)_t \cdot J_0))= \trace(A_{J_0X_t})(\phi_t^{-1}(x)),$$ as before in Equation \ref{trace with Upsilon}. Moreover, identically to Equation \ref{integral of trace equals via div to zero} we have \begin{equation}\label{integral over M of trace equals via div to zero}\int_M \trace(\Upsilon(x)_t (\Psi(x)_t \cdot J_0))\om^n(x) = 0.\end{equation}

 We shall now rewrite $\fS_{J_0}(\paph{\phi_t})$ via $\Psi_t(x)$. For all $x\in B$ we have the loop $\delta(x) = \paph{\Psi(x)_t \cdot J_0} \# \overline{[J_0, \Phi(x)_1 \cdot J_0]}$. We then for all $x \in B$ choose a disk $D(x)$ that bounds $\delta(x)$ - in fact one can construct $D(x)$ as the geodesic join of $\paph{\Psi(x)_t \cdot J_0}$ with $J_0$ - that is $D(x) = \bigcup_t \overline{[J_0, \Psi(x)_t \cdot J_0]}$ properly parametrized. Denote $\gamma_t(x)= \overline{[J_0, \Psi(x)_t \cdot J_0]}$. Denote by $\beta_t(x)$ the path $\{\Psi(x)_{t'} \cdot J_0\}_{t'=0}^{t'=t}$.

Compute \begin{equation}\label{write out fS via Psi begin} \fS_{J_0}(\paph{\phi_t}) = \int_D \Om - \intoi \mu(X_t)(\phi_t \cdot J_0) = \int_B (\int_{D(x)} \sigma_x) \om^n(x) - \intoi \int_M S(\phi_t \cdot J_0) H_t(x) \om^n(x).\end{equation} Now as before by Equation \ref{int_D sigma via tau_Lin and mu} and the definition of the moment map for the action of $G=Sp(2n)$ on $X=G/K$ \begin{equation}\label{local at x on M}\int_{D(x)} \sigma_x  \simeq \frac{1}{2}\tau_{Lin}(\paph{\Psi(x)_t})  - \int_0^1 f(x)_t(\Psi(x)_t \cdot J_0) dt\end{equation} where the function \begin{equation}\label{contact Hamiltonian for Psi} f(x)_t(J) = -\frac{1}{2} \trace(\Upsilon(x)_t J)\end{equation} is the contact Hamiltonian for the canonical lifting of $\Psi(x)_t$ to the principal $S^1$-bundle of unit vectors in $\Lambda^N_\C T\cS_x,\, N=n(n+1)/2$ (cf. Equation \ref{equation - moment map for finite dimensional groups}). Hence by Equations \ref{write out fS via Psi begin}, \ref{local at x on M}, \ref{contact Hamiltonian for Psi} and \ref{integral over M of trace equals via div to zero} we have \begin{equation}\label{write out fS via Psi} \fS_{J_0}(\paph{\phi_t}) \simeq  \frac{1}{2}\int_M \tau_{Lin}(\paph{\Psi(x)_t})\om^n(x) - \intoi \int_M S(\phi_t \cdot J_0) H_t(x) \om^n(x).\end{equation}

We shall now rewrite the Py quasimorphism $S_2$ from Definition \ref{Py quasimorphism} via $\Psi(x)_t$. Then comparing the effect of the difference in connections with the second term in Equation \ref{write out fS via Psi} we shall establish the equality.

First we note that the connection $\nabla$ gives us a parallel transport on $\cL(TM,\om)$ and on $P^2$, since it preserves $J_0$ and $\om$. Moreover, since the map $det^2:\cL(TM,\om) \to P^2$ is defined using only $J_0$ and $\om$ the following diagram commutes.

\begin{equation}\label{Lagrangian Grassmannian and parallel transport}
  \displaystyle    \begin{array}{clc}
         \mathcal{L}(TM,\om)_{(\phi_t)^{-1}x} & \xrightarrow{\Gamma_{\overline{p_{x,t}}}} & \mathcal{L}(TM,\om)_x\\
         \scriptstyle{det^2}{\downarrow} &   &  \scriptstyle{det^2}{\downarrow}\\
         P^2_{(\phi_t)^{-1}x} & \xrightarrow{\Gamma_{\overline{p_{x,t}}}} & P^2_x\\
      \end{array}
   \end{equation}
 In other words for $L_0 \in \cL_{(\phi_t)^{-1}x}(TM,\om)$ we have $det^2(\Gamma_{\overline{{p_{x,t}}}} L_0) = \Gamma_{\overline{{p_{x,t}}}} det^2(L_0)$. It will be more convenient to compute $S_2$ on the inverse path $\paph{\psi_t = \phi_t^{-1}}$. Indeed consider the paths $det^2({\psi_t}_{*_x} L)$ and $\widehat{\psi}_t (det^2(L))$ in $P^2$ for $L \in \cL(TM,\om)_x$. These paths differ by an angle as follows \[det^2({\psi_t}_{*_x}(L)) = e^{i2\pi \vartheta(t)} \widehat{\psi}_t(det^2(L))\] Then the paths $\Gamma_{\overline{p_{x,t}}} det^2({\psi_t}_{*_x} L) =  det^2(\Gamma_{\overline{p_{x,t}}} {\psi_t}_{*_x} L)$ (here we use Equation \ref{Lagrangian Grassmannian and parallel transport}) and $\Gamma_{\overline{p_{x,t}}} \widehat{\psi}_t (det^2(L))$ in $(P^2)_x$ also differ by the same angle. And since these are paths in one fiber, we have \[angle(L,\paph{\psi_t}) = varangle(\paph{e^{i2\pi\vartheta(t)}}) =\] \begin{equation}\label{Py angle L function via difference of varangles} varangle(\paph{det^2(\Gamma_{\overline{p_{x,t}}} {\psi_t}_{*_x} L)}) - varangle(\paph{\Gamma_{\overline{p_{x,t}}} \widehat{\psi}_t (det^2(L))}).\end{equation}

 Note that the second term in Equation \ref{Py angle L function via difference of varangles} does not depend on the choice of $L \in \cL(TM,\om)_x$, since both $\Gamma_{\overline{p_{x,t}}}$ and $\widehat{\psi}_t$ commute with rotations of the fibers. Therefore the function \[angle(x,\paph{\psi_t}) = \inf_{L \in \cL(TM,\om)_x} angle(L,\paph{\psi_t})\] satisfies \begin{equation}\label{Py angle x function via difference of varangles}angle(x,\paph{\psi_t}) = \inf_{L \in \cL(TM,\om)_x} (varangle(\paph{det^2(\Gamma_{\overline{p_{x,t}}} {\psi_t}_{*_x} L)})) - varangle(\paph{\Gamma_{\overline{p_{x,t}}} \widehat{\psi}_t y}),\end{equation} for any $y \in (P^2)_x$. Note first that ${\psi_t}_{*_x} = ({\phi_t}_{*_(\phi_t)^{-1}x})^{-1}$ and therefore $\Gamma_{\overline{p_{x,t}}} {\psi_t}_{*_x} = \Psi(x)_t^{-1}$. Then note that \begin{equation}\label{the Maslov part of Py's quasimorphism L}varangle(\paph{det^2(\Gamma_{\overline{p_{x,t}}} {\psi_t}_{*_x} L)}) \simeq \tau_{Lin} (\paph{\Psi(x)_t^{-1}}) = - \tau_{Lin}(\paph{\Psi(x)_t})\end{equation} by the construction of the Maslov quasimorphism on the universal cover of the linear symplectic group using its action on the Lagrangian Grassmannian \cite{BargeGhysEulerMaslov}. Therefore \begin{equation}\label{the Maslov part of Py's quasimorphism x}\inf_{L \in \cL(TM,\om)_x} (varangle(\paph{det^2(\Gamma_{\overline{p_{x,t}}} {\psi_t}_{*_x} L)})) \simeq - \tau_{Lin}(\paph{\Psi(x)_t}).\end{equation}
It remains to interpret the integral over $M$ with respect to $\om^n$ of the term $varangle(\paph{\Gamma_{\overline{p_{x,t}}} \widehat{\psi}_t y})$ in Equation \ref{Py angle x function via difference of varangles} via the Hermitian scalar curvature. For this purpose consider the two connection one-forms $\alpha$ and $\lambda$ on $P^2$ - where $d\alpha = 2 \til{\om}$ and $\lambda$ comes from the connection $\nabla$ on $TM$ and therefore satisfies $d\lambda = 2 \til{\rho}$ (for a form $\eta$ on $M$ we denote by $\til{\eta}$ its lift by the natural projection $P^2 \to M$). These connection one-forms differ by $\til{\theta} = \alpha - \lambda$ for a one-form $\theta$ on $M$. Then denoting by $Y_t$ the Hamiltonian vector field generating $\{\psi_t\}$ with normalized Hamiltonian $G_t$ (by the zero mean condition), and by $\wh{Y}_t$ the vector field generating $\{\wh{\psi}_t\}$ we have \begin{equation}\label{the Calabi part of Py's quasimorphism via thetas}varangle(\paph{\Gamma_{\overline{p_{x,t}}} \widehat{\psi}_t y}) = \intoi (\wh{\psi}_t)^*i_{\wh{Y}_t}\til{\theta}(x) dt = \intoi ({\psi}_t)^*i_{Y_t}\theta (x) dt.\end{equation}

We now compute as follows \begin{equation}\int_M \intoi ({\psi}_t)^*i_{Y_t}\theta(x) dt \om^n(x) = \intoi \int_M ({\psi}_t)^*i_{Y_t}\theta \om^n dt = \intoi \int_M i_{Y_t}\theta \om^n dt. \end{equation} It is therefore sufficient to compute the integrand \[\int_M i_{Y_t}\theta \om^n = n \int_M \theta i_{Y_t}\om \om^{n-1} = -n \int_M \theta d G_t \om^{n-1} = -n \int_M d\theta G_t \om^{n-1} =\] \[ = -2 n \int_M (\om - \rho) G_t \om^{n-1} = - 2n \int_M G_t \om^n + 2n \int_M G_t \rho \om^{n-1} =  2n \int_M G_t \rho \om^{n-1} = \] by the definition of the Hermitian scalar curvature \[ = 2\int_M G_t S(J_0) \om^n  = \] since $G_t(x) = -H_t(\phi_t x)$ by the cocycle formula \begin{equation}\label{the Calabi part of Py's quasimorphism via Hermitian scalar curvature}- 2\int_M S(J_0) H_t(\phi_t x) \om^n(x).\end{equation}

Therefore by Equations \ref{Py angle x function via difference of varangles}, \ref{the Maslov part of Py's quasimorphism x}, \ref{the Calabi part of Py's quasimorphism via thetas},\ref{the Calabi part of Py's quasimorphism via Hermitian scalar curvature} we have from the definition of $S_2$ (Definition \ref{Py quasimorphism}) that
\begin{equation}\label{S_2 of psi_t via Maslov and scalar curvature}-S_2(\paph{\psi_t}) \simeq -\int_M \tau_{Lin}(\paph{\Psi(x)_t}) \om^n(x) + 2\int_M S(J_0) H_t(\phi_t x) \om^n(x).\end{equation}

Therefore by Lemma \ref{nu of x and nu of x^{-1}} we have \begin{equation}\label{write out Py's quasimorphism via Psi} -S_2(\paph{\phi_t}) \simeq S_2(\paph{\psi_t}) = \int_M \tau_{Lin}(\paph{\Psi(x)_t}) \om^n(x) - 2\int_M S(J_0) H_t(\phi_t x) \om^n(x).\end{equation}

From Equations \ref{write out fS via Psi} and \ref{write out Py's quasimorphism via Psi} we conclude that \[2\fS_{J_0} \simeq -S_2,\] which by homogenizing gives \[2\fS = -\fS_{Py}\] finishing the proof.

\subsection{The restriction to the Entov quasimorphism}\label{Restriction to Entov}

Here we prove the second point of Theorem \ref{General agrees with Py and Entov} on the agreement of the general quasimorphism of Corollary \ref{Ham general principle} and the quasimorphism of Entov \cite{EntovCommutatorLength} from Definition \ref{Entov quasimorphism}. First we give an alternative definition of Entov's quasimorphism along the lines of the definition of Py's quasimorphism, which will more easily be shown to agree with the general quasimorphism.

\begin{df}\label{Entov quasimorphism 2}(A second definition of the quasimorphism $\fS_{En}$)
Given a symplectic manifold $(M,\om)$ with $c_1(TM,\om) = 0$ one first trivializes the top exterior power $\Lambda^n_\C (TM,J) \cong \C \times M$ of $(TM,\om,J)$ for $J \in \J$ as a Hermitian line bundle. The square $P^2$ of the unit frame bundle $S^1 \times M \cong P \xrightarrow{S^1} M$ of $L=\Lambda^n_\C(TM,J,\om)$ - that is the unitary frame bundle $P^2$ of $L^{\otimes 2}$ - admits a natural map $det^{2}: \cL(TM,\om) \to P^2$ from the Lagrangian Grassmannian bundle $\cL(TM,\om)$, since $\cL(TM)_x = U(TM_x,\om_x,J_x)/O(n)$. For a path $\arr{\phi}=\paph{\phi_t}$ in $\G$ with $\phi_0 = Id$, choosing a point $L \in \mathcal{L}(TM,\omega)_x$ we have the curve $\{{{\phi_t}_*}_x(L)\}_{0 \leq t \leq 1}$ in $\mathcal{L}(TM,\omega)$, and consequently the curve $\{ det^2({\phi_t}_{*_x}(L))\}_{0 \leq t \leq 1}$ in $P^2$. By means of the induced trivialization $P^2  \cong S^1 \times M$ this gives a continuous curve $e^{i2\pi \vartheta(t)}:[0,1] \to S^1$.  Define \[angle(L,\arr{\phi}) = varangle({\paph{e^{i2\pi \vartheta(t)}}}) = \vartheta(1) - \vartheta(0) ,\] and then the function \[angle(x,\arr{\phi}) = \displaystyle\inf_{L \in \mathcal{L}(TM,\omega)_x} angle(L,\arr{\phi})\] is measurable and bounded on $M$ and \[R_1(\arr{\phi}) = \int_M angle(x,\arr{\phi}) \om^n(x)\] does not depend on homotopies of $\arr{\phi}$ with fixed endpoints, defining a quasimorphism \[R_1: \til{\G} \to \R.\] Its homogeneization $\mathfrak{S}_{En}:\til{G} \to \R$, defined by $\mathfrak{S}_{En}(\til{\phi}):=\displaystyle\lim_{k \to \infty} \frac{R_1(\til{\phi}^k)}{k}$ is a homogenous quasimorhism on $\til{G}$ that is independent of the non-canonical choices of trivialization, and of the almost complex structure $J$.
\end{df}

\begin{prop}\label{Entov quasimorphism equivalence}
Definitions \ref{Entov quasimorphism} and \ref{Entov quasimorphism 2} for the Entov quasimorphism are equivalent.
\end{prop}

\begin{pfs}
Following Appendix C in \cite{LoopRemarks} one notes that the trivialization of $(TM,\om,J)$ over $U = M\setminus Z$ can be chosen to agree with the restriction from $M$ to $U$ of a given trivialization of $\Lambda^n_\C(TM,J)$. Then given a path $\arr{\phi}$ one immediately has $\simeq$ equality of the two $angle(x,\arr{\phi})$ functions on $U_{\arr{\phi}} = M \setminus Z_{\arr{\phi}}$ by the construction of the Maslov quasimorphism on the universal cover of the linear symplectic group using its action on the Lagrangian Grassmannian \cite{BargeGhysEulerMaslov} and the commutativity of the diagram

\begin{equation*}
      \begin{array}{clcccr}
         \mathcal{L}(TM,\om)|_U & \xrightarrow{det^2} & (\Lambda^n_\C(TM,J))^{\otimes 2}_U & \to & ^b\C|_U \\
         \scriptstyle{\wr}{\downarrow} &   &  \scriptstyle{\wr}{\downarrow} & & \shortparallel\\
         \mathcal{L}(\,^b\C^{n},\omega_{std})|_U & \xrightarrow{det^2} & (\bigwedge^{n}\,^b\C^{n})^{\otimes2}|_U & \cong & ^b\C|_U,
      \end{array}
   \end{equation*}

   where $^b\C$ is the trivial complex line bundle $\C \times M$ over $M$, and all vector bundles are complex and Hermitian.
\end{pfs}

Now we turn to showing the equality $\fS = \fS_{En}$. The proof is very similar to the one for the first point of Theorem \ref{General agrees with Py and Entov} and is even somewhat easier. Therefore we mostly outline the main steps and leave out details that are identical to those in Section \ref{Restriction to Py}.

First we recall Equation \ref{write out fS via Psi} \[\fS_{J_0}(\paph{\phi_t}) \simeq  \frac{1}{2}\int_M \tau_{Lin}(\paph{\Psi(x)_t})\om^n(x) - \intoi \int_M S(\phi_t \cdot J_0) H_t(x) \om^n(x).\]

We also recall the commutation relation of Equation \ref{Lagrangian Grassmannian and parallel transport}
\[
  \displaystyle    \begin{array}{clc}
         \mathcal{L}(TM,\om)_{(\phi_t)^{-1}x} & \xrightarrow{\Gamma_{\overline{p_{x,t}}}} & \mathcal{L}(TM,\om)_x\\
         \scriptstyle{det^2}{\downarrow} &   &  \scriptstyle{det^2}{\downarrow}\\
         P^2_{(\phi_t)^{-1}x} & \xrightarrow{\Gamma_{\overline{p_{x,t}}}} & P^2_x\\
      \end{array}
      \]
 That is for $L_0 \in \cL_{(\phi_t)^{-1}x}(TM,\om)$ we have $det^2(\Gamma_{\overline{{p_{x,t}}}} L_0) = \Gamma_{\overline{{p_{x,t}}}} det^2(L_0)$.

It will be more convenient to compute $R_1$ on the inverse path $\paph{\psi_t = \phi_t^{-1}}$. Indeed the path $det^2({\psi_t}_{*_x} L)$ in $P^2$ gives by the trivialization a smooth angle function $e^{i2\pi \vartheta(t)}:[0,1] \to S^1$. The path $\Gamma_{\overline{p_{x,t}}}: (P^2)_{(\phi_t)^{-1}x} \to (P^2)_x$ also gives by the trivialization a smooth angle function $e^{i2\pi \varphi(x,t)}: [0,1] \to S^1$. Noting the relation $\Gamma_{\overline{p_{x,t}}} det^2({\psi_t}_{*_x} L) =  det^2(\Gamma_{\overline{p_{x,t}}} {\psi_t}_{*_x} L)$ (by Equation \ref{Lagrangian Grassmannian and parallel transport}), we have \[angle(L,\paph{\psi_t}) = varangle(\paph{e^{i2\pi\vartheta(t)}}) =\] \begin{equation}\label{Entov angle L function via difference of varangles} varangle(\paph{det^2(\Gamma_{\overline{p_{x,t}}} {\psi_t}_{*_x} L)}) - varangle(\paph{e^{i2\pi \varphi(x,t)}}).\end{equation} Consequently, the function $angle(x,\paph{\psi_t}) = \inf_{L \in \cL(TM,\om)_x} angle(L,\paph{\psi_t})$ satisfies \begin{equation}\label{Entov angle x function via difference of varangles}angle(x,\paph{\psi_t}) = \inf_{L \in \cL(TM,\om)_x} (varangle(\paph{det^2(\Gamma_{\overline{p_{x,t}}} {\psi_t}_{*_x} L)})) - varangle(\paph{e^{i2\pi \varphi(x,t)}}).\end{equation} Note first that ${\psi_t}_{*_x} = ({\phi_t}_{*_(\phi_t)^{-1}x})^{-1}$ and therefore $\Gamma_{\overline{p_{x,t}}} {\psi_t}_{*_x} = \Psi(x)_t^{-1}$. Then note that \begin{equation}\label{the Maslov part of Entov's quasimorphism L}varangle(\paph{det^2(\Gamma_{\overline{p_{x,t}}} {\psi_t}_{*_x} L)}) \simeq \tau_{Lin} (\paph{\Psi(x)_t^{-1}}) = - \tau_{Lin}(\paph{\Psi(x)_t})\end{equation} by the construction of the Maslov quasimorphism on the universal cover of the linear symplectic group using its action on the Lagrangian Grassmannian \cite{BargeGhysEulerMaslov}. Therefore \begin{equation}\label{the Maslov part of Entov's quasimorphism x}\inf_{L \in \cL(TM,\om)_x} (varangle(\paph{det^2(\Gamma_{\overline{p_{x,t}}} {\psi_t}_{*_x} L)})) \simeq - \tau_{Lin}(\paph{\Psi(x)_t}).\end{equation}

It remains to interpret the integral over $M$ with respect to $\om^n$ of the term $varangle(\paph{e^{i2\pi \varphi(x,t)}})$ in Equation \ref{Entov angle x function via difference of varangles} via the Hermitian scalar curvature. For this purpose note that the trivialization $P^2 \cong S^1 \times M$ is equivalent to a flat connection $\alpha$ on $P^2$ without holonomy. Consider now the two connection one-forms $\alpha$ and $\lambda$ on $P^2$ - where in particular $d\alpha = 0$ and $\lambda$ comes from the connection $\nabla$ on $TM$ and therefore satisfies $d\lambda = 2 \til{\rho}$ (for a form $\eta$ on $M$ we denote by $\til{\eta}$ its lift by the natural projection $P^2 \to M$). These connection one-forms differ by $\til{\theta} = \alpha - \lambda$ for a one-form $\theta$ on $M$. Then denoting by $Y_t$ the Hamiltonian vector field generating $\{\psi_t\}$ with Hamiltonian $G_t$ normalized by the zero mean condition, and by $\wh{Y}_t$ the horizontal vector field that projects onto $Y_t$ generating the path $\{\wh{\psi}_t\}$ of $\alpha$-preserving diffeomorphism of $P^2$ (in other words $\wh{\psi_t} = Id \times \psi_t$ in the trivialization $P^2 \cong S^1 \times M$ ) we have \begin{equation}\label{the Calabi part of Entov's quasimorphism via thetas} varangle(\paph{e^{i2\pi \varphi(x,t)}}) = \intoi (\wh{\psi}_t)^*i_{\wh{Y}_t}\til{\theta}(x) dt = \intoi ({\psi}_t)^*i_{Y_t}\theta (x) dt.\end{equation}

We now compute as follows \begin{equation}\int_M \intoi ({\psi}_t)^*i_{Y_t}\theta(x) dt \om^n(x) = \intoi \int_M ({\psi}_t)^*i_{Y_t}\theta \om^n dt = \intoi \int_M i_{Y_t}\theta \om^n dt. \end{equation} It is therefore sufficient to compute the integrand \[\int_M i_{Y_t}\theta \om^n = n \int_M \theta i_{Y_t}\om \om^{n-1} = -n \int_M \theta d G_t \om^{n-1} = -n \int_M d\theta G_t \om^{n-1} =\] \[ = 2n \int_M \rho G_t \om^{n-1} = 2n \int_M G_t \rho \om^{n-1} =  2n \int_M G_t \rho \om^{n-1} = \] by the definition of the Hermitian scalar curvature \[ = 2\int_M G_t S(J_0) \om^n  = \] since $G_t(x) = -H_t(\phi_t x)$ by the cocycle formula \begin{equation}\label{the Calabi part of Entov's quasimorphism via Hermitian scalar curvature}= - 2\int_M S(J_0) H_t(\phi_t x) \om^n(x).\end{equation}

Therefore by Equations \ref{Entov angle x function via difference of varangles}, \ref{the Maslov part of Entov's quasimorphism x}, \ref{the Calabi part of Entov's quasimorphism via thetas},\ref{the Calabi part of Entov's quasimorphism via Hermitian scalar curvature} we have from the definition of $R_1$ that
\begin{equation}\label{R_1 of psi_t via Maslov and scalar curvature}R_1(\paph{\psi_t}) \simeq -\int_M \tau_{Lin}(\paph{\Psi(x)_t}) \om^n(x) + 2\int_M S(J_0) H_t(\phi_t x) \om^n(x).\end{equation}

Therefore by Lemma \ref{nu of x and nu of x^{-1}} we have \begin{equation}\label{write out Entov's quasimorphism via Psi} R_1(\paph{\phi_t}) \simeq -R_1(\paph{\psi_t}) = \int_M \tau_{Lin}(\paph{\Psi(x)_t}) \om^n(x) - 2\int_M S(J_0) H_t(\phi_t x) \om^n(x).\end{equation}

From Equations \ref{write out fS via Psi} and \ref{write out Entov's quasimorphism via Psi} we conclude that \[2\fS_{J_0} \simeq R_1,\] which by homogenizing gives \[2\fS = \fS_{En}\] finishing the proof.

\subsection{Calibrating the $L^2_2$ norm}\label{calibrating the L^2_2 norm}

Here we derive Equation \ref{bound by L^2_2 of fS}.

Note that the second summand of $\fS_{J_0}(\arr{\phi}) = \int_{D_{\arr{\phi}}} \Om - \intoi  S(J_0) H_t(\phi_t x) \om^n(x) dt$ satisfies \begin{equation}\label{bound by L^p_k of the second part of fS}|\intoi S(J_0) H_t(\phi_t x) \om^n(x) dt| \leq ||S(J_0)||_{L^q(M,\om^n)} \cdot \intoi ||H_t||_{L^p(M,\om^n)} dt\end{equation} where $1 \leq p,q \leq \infty$ and $1/p+1/q = 1$ and is therefore bounded by $C_{p} ||\arr{\phi_t}||_{k,p}$ for every $k \geq 0$ and $1 \leq p \leq \infty$.

Let us turn to the first summand $\int_{D_{\arr{\phi}}} \Om$. First note that since on the Siegel upper half-space $\cS_n$ the natural invariant Kahler form $\sigma_{\cS_n}$ has a primitive $\lambda_{\cS_n}$ that is bounded by a constant $C(n)$ with respect to the metric induced by $(\sigma_{\cS_n},j_{\cS_n})$ and vanishes on geodesics starting at $iId$, the infinite-dimensional space $(\J,\Om,\mathbb{J})$ also has a primitive $\Lambda$ for $\Om$ that is bounded with respect to the metric induced by $(\Om,\mathbb{J})$ by the constant $C(n,\om) = C(n) \mVol^{1/2}$ and vanishes on geodesics starting at $J_0$. That is \[|\Lambda(\Upsilon)| \leq C(n) \mVol^{1/2} \Om(\Upsilon,\mathbb{J}\Upsilon)^{1/2},\] for a vector $\Upsilon \in T_J \J$. In that case $\int_{D_{\arr{\phi}}} \Om = \int_{\arr{\phi}\cdot J_0} \Lambda = \intoi \Lambda_{\phi_t \cdot J_0} ((\phi_t)_* L_{X_t} J_0) dt$ and consequently \[|\int_{D_{\arr{\phi}}} \Om| \leq C(n,\om) \intoi \Om_{\phi_t \cdot J_0}((\phi_t)_* L_{X_t}J_0,(\phi_t \cdot J_0) ((\phi_t)_* L_{X_t}J_0))^{1/2}dt \leq\]  \[\leq C'(n,\om) \intoi (\int_M trace((\phi_t)_* (L_{X_t}J_0)^2) \om^n)^{1/2} dt=\] \[ = C'(n,\om) \intoi (\int_M trace((L_{X_t}J_0)^2) \om^n)^{1/2}dt \leq \] \begin{equation}\label{bound by L^2_2 of the first part of fS}\leq C''(n,\om,J_0) \intoi (\int_M (|X_t|^2 + |\nabla X_t|^2) \om^n)^{1/2} dt\leq C^{(3)}(n,\om,J_0) ||\arr{\phi}||_{2,2}.\end{equation}

Therefore by Equations \ref{bound by L^2_2 of the first part of fS} and \ref{bound by L^p_k of the second part of fS} we have for all $\til{\phi} \in \til{\G}$ the estimate \begin{equation}\fS_{J_0}(\til{\phi}) \leq C(n,\om,J_0) ||\til{\phi}||_{2,2}.\end{equation}

\section{Discussion}

\begin{enumerate}

\item  It was shown by Donaldson in \cite{DonaldsonMomentMapsDiffeomorphisms,DonaldsonMomentMapsDiffGeom} that $\G$ acts in a Hamiltonian way on additional spaces (e.g. spaces of submanifolds/cycles). These may yield more homomorphisms $\pi_1(Ham) \to \R$ by the Action-homomorphism construction for equivariant moment maps, and perhaps new quasimorphisms on $\til{\G}$. Moreover, Futaki shows in \cite{FutakiMomentMaps} that the space $\J_{int} \subset \J$ of integrable almost complex structures can be endowed with additional symplectic structures that give different moment maps for the action of $\G$, from which the Bando-Futaki invariants $F_{c_k}$ are obtained when restricting to the subgroup $\G_{J_0}$. It would be interesting to extend the methods of Futaki to all $\J$, taking care of the Nijenhuis tensor, and to check two things. First it is most likely that the corresponding Action-homomorphisms on $\pi_1(\G)$ will coincide with the invariants $I_{c_k}$ (cf. \cite{symplectures}) obtained by integrating the $k$-th vertical Chern class times $u^{n+1-k}$ in Definition \ref{generalized Action-Maslov c_1}. Second, it would be interesting to extend the perturbation of Futaki to incorporate such invariants as $I_{c_1 c_2}$ corresponding to symmetric polynomials that are not elementary.

\item It is interesting to note that the Entov quasimorphism (Defintions \ref{Entov quasimorphism}, \ref{Entov quasimorphism 2}) is defined on the extension $\cH = Symp(M,\om)$ of the group $\G = Ham(M,\om)$, while the moment map picture is currently stated for the action of $\G$ on $\J$ only. It is therefore interesting to check whether in the case $c_1(TM,\om) = 0$ the moment map for the action of $\G$ on $\J$ extends to a moment map for the action of $\cH$ on $\J$ - along the lines of \cite{DonaldsonMomentMapsDiffGeom} for example. It would also be interesting to investigate the possibility of extending the moment map this way without conditions on $c_1(TM,\om)$ - to provide an extension when it is possible and to investigate the obstructions to extending when the extension is not possible. This may well be related to the Flux homomorphism.
\item It is interesting to investigate the restriction $\A_\mu$ of $\fS$ to $\pi_1 \G$ for symplectic manifolds $(M,\om)$ of finite volume that are not closed. Does this restriction have an interpretation like $I_{c_1}$ in terms of characteristic numbers of the associated Hamiltonian vector bundle? It would also be interesting to say something new about the Entov quasimorphism in the new interpretation - can it be computed for example for the new symplectic manifolds constructed by Fine and Panov (\cite{PanovFineSurvey} and references therein)?
\item It would be interesting to compare the general principle for generating quasimorphisms introduced in this paper with other general constructions of quasimorphisms. While the relation to the Burger-Iozzi-Wienhard construction of the rotation number from \cite{BurgerIozziWienhardMaximalToledoInvariant} is at least intuitively relatively simple to trace, the relation to the works of Ben Simon and Hartnick \cite{GabiTobias1,GabiTobias2,GabiTobias3} (cf. Calegari \cite{CalegariScl}) is somewhat more mysterious, since there seems to be no straightforward analogue of the Shilov boundary for the space $(\J,\Om,\mathbb{J})$ of compatible almost complex structures on $(M,\om)$. Hence it is an interesting question to exhibit a specific explicit invariant partial order or poset that gives the quasimorphism $\mathfrak{S}$ on $\til{Ham}(M,\om)$.
\item From a general philosophical point of view the action of $\G=Ham(M,\om)$ on $\J$ with Donaldson's equivariant moment map allows one to consider $\G$ in its $C^1$-topology as a generalized Hermitian Lie group with a generalized Hermitian symmetric space of non-compact type - in a way it behaves similarly to $Sp(2n,\R)$, which would be a "Hermitian" feature of $\G$. In comparison, the group $\G$ with the Hofer metric and related invariants is known to exhibit certain "hyperbolic features" (cf. \cite{PolterovichDynamicsGroups}) - shared with Gromov-hyperbolic finitely generated groups. This approach can be used to study the representations into $\G$ of fundamental groups of compact Kahler manifolds - e.g. Riemann surfaces of genus at least $2$. It is easy to construct the analogue of the Toledo invariant for representations of surface groups (using the bounded 2-cocycle of Reznikov \cite{ReznikovCharClassSymp,ReznikovCocyclesVol,ReznikovAnalyticTop} that equals the differential of $\mathfrak{S}_{J_0}$ - which corresponds to the "bounded Kahler class") that satisfies a corresponding Milnor-Wood type inequality (this can for example be proven using the quasimorphism $\fS_{J_0}$). One could then check which values of the Toledo invariant can be attained - note that this value will be $I_{c_1}$ on a certain loop $\gamma_\rho$ associated with the representation $\rho$, and hence for Kahler-Einstein manifolds is conjectured to vanish \cite{LoopRemarks} - this holds for example on $(\C P^n,\om_{FS})$ \cite{EntovPolterovichCalabiQmQh,RigidSubsets}. These methods could possibly be used to obtain restrictions on Hamiltonian actions of such groups, which would be complementary to those established by Polterovich (cf. \cite{PolterovichDynamicsGroups}), since surface groups are undistorted. In particular the notion of maximal representations (following works of Burger-Iozzi-Wienhard and others cf. \cite{BurgerIozziWienhardSurvey} for a survey) could be defined and their properties studied. The above-mentioned works of Ben Simon and Hartnick could again be of some use.

    Note also that while certain embeddings of right-angled Artin groups (and hence of most surface groups) into $\G$  of any symplectic manifold were constructed by Kapovich in \cite{KapovichRAAGs} these representations will have zero Toledo invariant. Indeed these constructions either factor through the subgroup $\G_B$ of diffeomorphisms supported in a ball, where the restriction of the quasimorphism to $\pi_1$ is trivial (c.f. definition \ref{Barge-Ghys average Maslov quasimorphism} of the Barge-Ghys average Maslov quasimorphism) or take values in $\G$ of a surface of genus $g$, where the restriction $I_{c_1}$ vanishes since $\pi_1(\G)$ is trivial (or torsion for the sphere). The surface can also have boundary - the Toledo invariant will still vanish, by the embedding functoriality (Proposition \ref{embedding functoriality}). However, it is quite an easy fact that since $\til{Ham}(M,\om)$ for closed $M$ is perfect by a theorem of Banyaga \cite{BanyagaPerfectSimple}, every element $\gamma \in \pi_1 Ham$ is of the form $\gamma=\gamma_\rho$ for some representation $\rho: \pi_1 (\Sigma_g) \to Ham$ (one can take $g$ to be the commutator length of $\gamma \in \til{Ham}$). Hence for $M = Bl_1(\C P^1)$, say, there is a nonzero Toledo invariant representation, the corresponding class in $\pi_1$ represented by a toric loop. It would therefore be interesting to write this class explicitly as a product of commutators in $\til{Ham}$.

\item Another interesting computation to make is that of $\fS$ on Hamiltonian paths generated by a time-independent (autonomous) Hamiltonian. This would give a \emph{quasi-state}-type functional (cf. e.g. \cite{EntovPoletrovichSymplecticQuasistatesSemiSimplicity,EntovPoletrovichLieQuasiStates}) on $C^\infty({M},{\R})$ corresponding to the quasimorphism $\fS$. This functional would retain the properties of linearity on Poisson-commutative subspaces and $Symp(M,\om)$-invariance, however it would not be monotone (since this would imply continuity in the $L^\infty$-norm) or vanish on functions with supports displaceable by Hamiltonian isotopies. In particular, it would be curious to find a formula for the value of this quasi-state on Morse functions on the manifold in terms of local data around the critical points, similarly to what was computed by Py in his thesis \cite{PyThesis} for the case of the two-sphere $S^2$. Here Equation \ref{write out fS via Psi} could be very useful. One could also ask whether there are similar localization formulas for actions of other groups e.g. $\R^k$ with tame fixed manifolds. For one, in the case when $(M,\om)$ is toric the restriction $I_{c_1}$ of $\fS$ to $\pi_1(\G)$ has been computed on loops coming from the torus action (cf. \cite{LoopRemarks} and references therein).
\end{enumerate}

\nocite{DonaldsonBoundsOnCalabi}

\bibliographystyle{amsplain}
\bibliography{ActionFunctionalAndMomentMapsRefs}

\providecommand{\bysame}{\leavevmode\hbox to3em{\hrulefill}\thinspace}
\providecommand{\MR}{\relax\ifhmode\unskip\space\fi MR }
\providecommand{\MRhref}[2]{%
  \href{http://www.ams.org/mathscinet-getitem?mr=#1}{#2}
}
\providecommand{\href}[2]{#2}
\begin{thebibliography}{10}

\bibitem{AbreuGranjaKitchlooSketch}
M.~Abreu, G.~Granja, and N.~Kitchloo, \emph{Moment maps, symplectomorphism
  groups and compatible complex structures}, J. Symplectic Geom. \textbf{3}
  (2005), no.~4, 655–--670, Special Issue. Conference on Symplectic Topology.

\bibitem{AbreuGranjaKitchlooBIG}
\bysame, \emph{Compatible complex structures on symplectic rational ruled
  surfaces}, Duke Math. J. \textbf{148} (2009), no.~3, 539--–600.

\bibitem{AtiyahBottConnectionsPrelim}
M.~F. Atiyah and R.~Bott, \emph{{Y}ang-{M}ills and bundles over algebraic
  curves}, Geometry and analysis, vol.~5, Indian Acad. Sci., Bangalore, 1980,
  pp.~11--20.

\bibitem{AtiyahBottConnections}
\bysame, \emph{The {Y}ang-{M}ills equations over {R}iemann surfaces}, Philos.
  Trans. Roy. Soc. London Ser. A \textbf{308} (1983), no.~1505, 523--–615.

\bibitem{BanyagaPerfectSimple}
A.~Banyaga, \emph{Sur la structure du groupe des
  diff$\acute{\text{e}}$omorphismes qui pr$\acute{\text{e}}$servent une forme
  symplectique}, Comment. Math. Helv. \textbf{53} (1978), no.~2, 174--–227.

\bibitem{BargeGhysEulerBornee}
J.~Barge and $\acute{\text{E}}$. Ghys, \emph{Surfaces et cohomologie
  born$\acute{\text{e}}$e}, Invent. Math. \textbf{92} (1988), no.~3, 509–--526.

\bibitem{BargeGhysEulerMaslov}
\bysame, \emph{Cocycles d'{E}uler et de {M}aslov}, Math. Ann. \textbf{294}
  (1992), 235--265.

\bibitem{GabiCalabiQm}
G.~{Ben Simon}, \emph{The nonlinear {M}aslov index and the {C}alabi
  homomorphism}, Commun. Contemp. Math. \textbf{9} (2007), no.~6, 769--780.

\bibitem{BenaimGambaudoMetricProperties}
M.~Benaim and J.-M. Gambaudo, \emph{Metric properties of the group of area
  preserving diffeomorphisms}, Trans. Amer. Math. Soc. \textbf{353} (2001),
  no.~11, 4661--–4672.

\bibitem{BottGr}
R.~Bott, \emph{The geometry and representation theory of compact {L}ie groups
  (notes by {G}.{L}. {L}uke)}, Representation theory of Lie groups. Proceedings
  of the SRC/LMS Research Symposium held in Oxford, June 28--July 15, 1977.
  (Cambridge-New York) (G.~L. Luke, ed.), London Mathematical Society Lecture
  Note Series, vol.~34, Cambridge University Press, 1979, pp.~65--90.

\bibitem{BrandThesis}
M.~Brandenbursky, \emph{Knot invariants and their applications to constructions
  of quasi-morphisms on groups}, Ph.D. thesis, Technion - Israel Institute of
  Technology, 2010, http://www.math.vanderbilt.edu/~brandem/Thesis-MB.pdf.

\bibitem{BrandenburskyEstimates}
\bysame, \emph{Quasi-morphisms and ${L}^p$-metrics on groups of
  volume-preserving diffeomorphisms}, preprint arXiv:1110.3353, 2011.

\bibitem{BurgerIozziWienhardSurvey}
M.~Burger, A.~Iozzi, and A.~Wienhard, \emph{Higher
  {T}eichm$\ddot{\text{u}}$ller spaces: from ${SL}(2,{\R})$ to other {L}ie
  groups}, preprint arXiv:1004.2894v2, 2010.

\bibitem{BurgerIozziWienhardMaximalToledoInvariant}
\bysame, \emph{Surface group representations with maximal {T}oledo invariant},
  Ann. of Math. (2) \textbf{172} (2010), no.~1, 517--–566.

\bibitem{BurgerMonodLattices}
M.~Burger and N.~Monod, \emph{Bounded cohomology of lattices in higher rank
  {L}ie groups}, J. Eur. Math. Soc. (JEMS) \textbf{1} (1999), no.~2, 199--–235.

\bibitem{CalabiCY}
E.~Calabi, \emph{On {K}\"{a}hler manifolds with vanishing canonical class.},
  Algebraic geometry and topology. A symposium in honor of S. Lefschetz
  (Princeton, N. J.), Princeton University Press.

\bibitem{CalabiHomomorphism}
\bysame, \emph{On the group of automorphisms of a symplectic manifold},
  Problems in analysis (Lectures at the Sympos. in honor of Salomon Bochner,
  Princeton Univ., Princeton, N.J., 1969), Princeton Univ. Press, 1970,
  pp.~1--26.

\bibitem{CalabiExtremalMetrics2}
\bysame, \emph{Extremal {K}$\ddot{\text{a}}$hler metrics. {I}{I}}, Seminar on
  Differential Geometry, Ann. of Math. Stud., vol. 102, Princeton Univ. Press,
  Princeton, N.J., 1982, pp.~259–--290.

\bibitem{CalabiExtremalMetrics1}
\bysame, \emph{Extremal {K}$\ddot{\text{a}}$hler metrics.}, Differential
  geometry and complex analysis, Springer, Berlin, 1985, pp.~95–--114.

\bibitem{CalegariScl}
D.~Calegari, \emph{scl}, MSJ Memoirs, vol.~20, Mathematical Society of Japan,
  Tokyo, 2009.

\bibitem{EilenbergChevalley}
C.~Chevalley and S.~Eilenberg, \emph{Cohomology theory of {L}ie groups and
  {L}ie algebras}, Trans. Amer. Math. Soc. \textbf{63} (1948), 85--124.

\bibitem{ClercKoufany}
J.-L. Clerc and K.~Koufany, \emph{Primitive du cocycle de {M}aslov
  g$\acute{\text{e}}$n$\acute{\text{e}}$ralis$\acute{\text{e}}$}, Math. Ann.
  \textbf{337} (2007), no.~1, 91–--138.

\bibitem{ClercOrsted}
J.-L. Clerc and B.~${\O}$rsted, \emph{The {G}romov norm of the {K}aehler class
  and the {M}aslov index}, Asian J. Math. \textbf{7} (2003), no.~2, 269–--295.

\bibitem{DomicToledo}
A.~Domic and D.~Toledo, \emph{The {G}romov norm of the {K}aehler class of
  symmetric domains}, Math. Ann. \textbf{276} (1987), no.~3, 425--–432.

\bibitem{DonaldsonComplexSurfacesHitchinKobayashi}
S.~K. Donaldson, \emph{Anti self-dual {Y}ang-{M}ills connections over complex
  algebraic surfaces and stable vector bundles}, Proc. London Math. Soc. (3)
  \textbf{50} (1985), no.~1, 1--26.

\bibitem{DonaldsonInfiniteDeterminantsHitchinKobayashi}
\bysame, \emph{Infinite determinants, stable bundles and curvature}, Duke Math.
  J. \textbf{54} (1987), no.~1, 231–--247.

\bibitem{DonaldsonGaugeTheoryComplexGeometry}
\bysame, \emph{Remarks on gauge theory, complex geometry and $4$-manifold
  topology}, Fields Medallists' lectures, World Sci. Ser. 20th Century Math.,
  vol.~5, World Sci. Publ., River Edge, NJ, 1997, pp.~384--403.

\bibitem{DonaldsonMomentMapsDiffeomorphisms}
\bysame, \emph{Moment maps and diffeomorphisms}, Asian J. Math. \textbf{3}
  (1999), no.~1, 1--15, Special Issue. Sir Michael Atiyah: a great
  mathematician of the twentieth century.

\bibitem{DonaldsonMomentMapsDiffGeom}
\bysame, \emph{Moment maps in differential geometry}, Surveys in differential
  geometry Vol. VIII (Boston, MA, 2002), Surv. Differ. Geom., vol. VIII, Int.
  Press, Somerville, MA, 2003, pp.~171–--189.

\bibitem{DonaldsonBoundsOnCalabi}
\bysame, \emph{Lower bounds on the {C}alabi functional}, J. Differential Geom.
  \textbf{70} (2005), no.~3, 453–--472.

\bibitem{DupontGuichardet}
J.~L. Dupont and A.~Guichardet, \emph{{$\grave{\text{A}}$} propos de l'article:
  "{S}ur la cohomologie r$\acute{\text{e}}$elle des groupes de {L}ie simples
  r$\acute{\text{e}}$els'' par {G}uichardet et {D}. {W}igner.}, Ann. Sci.
  $\acute{\text{E}}$cole Norm. Sup. (4) \textbf{11} (1978), no.~2, 293--–295.

\bibitem{EliashbergRatiu}
Y.~Eliashberg and T.~Ratiu, \emph{The diameter of the symplectomorphism group
  is infinite}, Invent. Math. \textbf{103} (1991), no.~2, 327–--340.

\bibitem{EntovCommutatorLength}
M.~Entov, \emph{Commutator length of symplectomorphisms}, Comment. Math. Helv.
  \textbf{79} (2004), no.~1, 58--104.

\bibitem{EntovPolterovichCalabiQmQh}
M.~Entov and L.~Polterovich, \emph{Calabi quasimorphism and quantum homology},
  Int. Math. Res. Not. (2003), no.~30, 1635--1676.

\bibitem{EntovPoletrovichSymplecticQuasistatesSemiSimplicity}
\bysame, \emph{Symplectic quasi-states and semi-simplicity of quantum
  homology}, Toric Topology (M.Masuda M.Harada, Y.Karshon and T.Panov, eds.),
  Contemporary Mathematics, vol. 460, Amer. Math. Soc., 2008,
  arXiv:0705.3735v3, pp.~47--70.

\bibitem{RigidSubsets}
\bysame, \emph{Rigid subsets of symplectic manifolds}, Compos. Math.
  \textbf{145} (2009), no.~3, 773--826.

\bibitem{EpsteinFujiwara}
D.~B. Epstein and K.~Fujiwara, \emph{The second bounded cohomology of
  word-hyperbolic groups.}, Topology \textbf{36} (1997), no.~6, 1275–--1289.

\bibitem{Fine}
J.~Fine, \emph{The {H}amiltonian geometry of the space of unitary connections
  with symplectic curvature}, preprint arXiv:1101.2420v1, 2011.

\bibitem{PanovFineSurvey}
J.~Fine and D.~Panov, \emph{Building symplectic manifolds using hyperbolic
  geometry}, Proceedings of the G$\ddot{\text{o}}$kova Geometry-Topology
  Conference 2009, Int. Press, Somerville, MA, 2010, pp.~124--136.

\bibitem{FomenkoFuks}
A.~T. Fomenko and D.~B. Fuks, \emph{Kurs gomotopicheskoi topologii. ({R}ussian)
  [{A} course in homotopic topology]}, "Nauka", Moscow, 1989.

\bibitem{Fujiki}
A.~Fujiki, \emph{Moduli space of polarized algebraic manifolds and
  {K}$\ddot{\text{a}}$hler metrics}, Sugaku Expositions \textbf{5} (1992),
  no.~2, 173–--191.

\bibitem{Fujiwara}
K.~Fujiwara, \emph{The second bounded cohomology of a group acting on a
  {G}romov-hyperbolic space}, Proc. London Math. Soc. (3) \textbf{76} (1998),
  no.~1, 70--–94.

\bibitem{f}
A.~Futaki, \emph{An obstruction to the existence of {E}instein {K}\"{a}hler
  metrics}, Invent. Math. \textbf{73} (1983), no.~3, 437--443.

\bibitem{FutakiMomentMaps}
\bysame, \emph{Harmonic total {C}hern forms and stability}, Kodai Math. J.
  \textbf{29} (2006), no.~3, 346–--369.

\bibitem{GhysGambaudoEnlacements}
J.-M. Gambaudo and $\acute{\text{E}}$. Ghys, \emph{Enlacements asymptotiques},
  Topology \textbf{36} (1997), no.~6, 1355–--1379.

\bibitem{GhysGambaudoCommutatorsSurfaces}
\bysame, \emph{Commutators and diffeomorphisms of surfaces}, Ergodic Theory
  Dynam. Systems \textbf{24} (2004), no.~5, 1591–--1617.

\bibitem{GambaudoLagrange}
J.-M. Gambaudo and M.~Lagrange, \emph{Topological lower bounds on the distance
  between area preserving diffeomorphisms}, Bol. Soc. Brasil. Mat. (N.S.)
  \textbf{31} (2000), no.~1, 9–--27.

\bibitem{Gauduchon}
P.~Gauduchon, \emph{Hermitian connections and {D}irac operators}, Boll. Un.
  Mat. Ital. B (7) \textbf{11} (1997), no.~2, suppl., 257–--288.

\bibitem{GiventalQuasimorphism}
A.~Givental, \emph{Nonlinear generalization of the {M}aslov index}, Theory of
  singularities and its applications, Adv. Soviet Math., vol.~1, Amer. Math.
  Soc., Providence, RI, 1990, pp.~71--103.

\bibitem{GromovPseudohol}
M.~Gromov, \emph{Pseudoholomorphic curves in symplectic manifolds}, Invent.
  Math. \textbf{82} (1985), no.~2, 307--347.

\bibitem{GromovSandH}
\bysame, \emph{Soft and hard symplectic geometry}, Proceedings of the
  {I}nternational {C}ongress of {M}athematicians, {V}ol. 1, 2 ({B}erkeley,
  {C}alif., 1986) (Providence, RI), Amer. Math. Soc., 1987, pp.~81–--98.

\bibitem{GuichardetWigner}
A.~Guichardet and D.~Wigner, \emph{Sur la cohomologie r$\acute{\text{e}}$elle
  des groupes de {L}ie simples r$\acute{\text{e}}$els}, Ann. Sci.
  $\acute{\text{E}}$cole Norm. Sup. (4) \textbf{11} (1978), no.~2, 277--–292.

\bibitem{HamenstadtDiscrete}
U.~Hamenst$\ddot{\text{a}}$dt, \emph{Bounded cohomology and isometry groups of
  hyperbolic spaces}, J. Eur. Math. Soc. (JEMS) \textbf{10} (2008), no.~2,
  315–--349.

\bibitem{Hamenstadt}
\bysame, \emph{Isometry groups of proper hyperbolic spaces}, Geom. Funct. Anal.
  \textbf{19} (2009), no.~1, 170–--205.

\bibitem{Helgason}
S.~Helgason, \emph{Differential geometry, {L}ie groups, and symmetric spaces},
  Pure and Applied Mathematics, vol.~80, Academic Press, Inc., New York-London,
  1978.

\bibitem{KapovichRAAGs}
M.~Kapovich, \emph{R{A}{A}{G}s in {H}am}, preprint arXiv:1104.0348v1, 2011.

\bibitem{Knapp}
A.~W. Knapp, \emph{Lie groups beyond an introduction}, second ed., Progress in
  Mathematics, vol. 140, Birkh$\ddot{\text{a}}$user Boston, Inc., Boston, MA,
  2002.

\bibitem{KobayashiConnectionsAlmsotComplex}
S.~Kobayashi, \emph{Natural connections in almost complex manifolds},
  Explorations in complex and {R}iemannian geometry, Contemp. Math., vol. 332,
  Amer. Math. Soc., Providence, RI, 2003, pp.~153--–169.

\bibitem{KobayashiNomizu}
S.~Kobayashi and K.~Nomizu, \emph{{F}oundations of differential geometry {V}ol.
  {I}.}, Interscience Publishers, a division of John Wiley \& Sons, New
  York-London, 1963.

\bibitem{LMP}
F.~Lalonde, D.~McDuff, and L.~Polterovich, \emph{Topological rigidity of
  {H}amiltonian loops and quantum homology}, Invent. Math. \textbf{135} (1999),
  no.~2, 369--385.

\bibitem{symplectures}
D.~McDuff, \emph{Lectures on groups of symplectomorphisms}, Rend. Circ. Mat.
  Palermo (2) Suppl. (2004), no.~72, 43--78.

\bibitem{HamTracesContract}
\bysame, \emph{A survey of the topological properties of symplectomorphism
  groups.}, Topology, geometry and quantum field theory (Cambridge), London
  Math. Soc. Lecture Note Ser., vol. 308, Cambridge University Press, 2004,
  pp.~173--193.

\bibitem{IntroSymp}
D.~McDuff and D.~Salamon, \emph{Introduction to symplectic topology}, second
  ed., Oxford Mathematical Monographs, The Clarendon Press, Oxford University
  Press, New York, 1998.

\bibitem{JHolSymp}
\bysame, \emph{{$J$}-holomorphic curves and symplectic topology}, American
  Mathematical Society Colloquium Publications, vol.~52, American Mathematical
  Society, Providence, RI, 2004.

\bibitem{MilnorInfiniteDimensional}
J.~Milnor, \emph{Remarks on infinite-dimensional {L}ie groups}, Relativity,
  groups and topology, II (Les Houches, 1983), North-Holland, Amsterdam, 1984,
  p.~1007–1057.

\bibitem{Mok}
N.~Mok, \emph{Metric rigidity theorems on {H}ermitian locally symmetric
  manifolds}, Series in Pure Mathematics, vol.~6, World Scientific Publishing
  Co., Inc., Teaneck, NJ, 1989.

\bibitem{MonodShalom}
N.~Monod and Y.~Shalom, \emph{Cocycle superrigidity and bounded cohomology for
  negatively curved spaces}, J. Differential Geom. \textbf{67} (2004), no.~3,
  395--–455.

\bibitem{Omori}
H.~Omori, \emph{Infinite-dimensional {L}ie groups}, Translations of
  Mathematical Monographs, vol. 158, American Mathematical Society, Providence,
  RI, 1997.

\bibitem{YaronCalabiQmNonnMonotone}
Y.~Ostrover, \emph{Calabi quasi-morphisms for some non-monotone symplectic
  manifolds}, Algebr. Geom. Topol. \textbf{6} (2006), 405--–434.

\bibitem{YaronTyomkin}
Y.~Ostrover and I.~Tyomkin, \emph{On the quantum homology algebra of toric
  {F}ano manifolds}, Selecta Math. (N.S.) \textbf{15} (2009), no.~1, 121--149.

\bibitem{PolterovichLoops}
L.~Polterovich, \emph{Hamiltonian loops and {A}rnold's principle.}, Topics in
  singularity theory (Providence, RI), vol. 180, Amer. Math. Soc., 1997,
  pp.~181--187.

\bibitem{PolterovichGeometryGroupSymp}
\bysame, \emph{The geometry of the group of symplectic diffeomorphisms},
  Lectures in Mathematics ETH Z$\ddot{\text{u}}$rich,
  Birkha$\ddot{\text{u}}$ser Verlag, Basel, 2001.

\bibitem{PolterovichDynamicsGroups}
\bysame, \emph{Floer homology, dynamics and groups}, Morse theoretic methods in
  nonlinear analysis and in symplectic topology (Dordrecht), NATO Sci. Ser. II
  Math. Phys. Chem., vol. 217, Springer, 2006, pp.~417–--438.

\bibitem{EntovPoletrovichLieQuasiStates}
L.~Polterovich and M.~Entov, \emph{Lie quasi-states}, Journal of Lie Theory
  \textbf{19} (2009), no.~3, 613--637.

\bibitem{pyqm}
P.~Py, \emph{Quasi-morphismes et invariant de {C}alabi}, Ann. Sci. E'cole Norm.
  Sup. (4) \textbf{39} (2006), no.~1, 177--195.

\bibitem{PyThesis}
\bysame, \emph{Quasi-morphismes et diff$\grave{\text{e}}$omorphismes
  {H}amiltoniens}, Ph.D. thesis, $\grave{\text{E}}$cole normale superieure de
  Lyon, France, 2008, http://tel.archives-ouvertes.fr/tel-00263607/fr/.

\bibitem{ReznikovCharClassSymp}
A.~Reznikov, \emph{Characteristic classes in symplectic topology. {A}ppendix
  {D} by {L}udmil {K}atzarkov.}, Selecta Math. (N.S.) \textbf{3} (1997), no.~4,
  601–--642.

\bibitem{ReznikovCocyclesVol}
\bysame, \emph{Continuous cohomology of the group of volume-preserving and
  symplectic diffeomorphisms, measurable transfer and higher asymptotic
  cycles}, Selecta Math. (N.S.) \textbf{5} (1999), no.~1, 181--–198.

\bibitem{ReznikovAnalyticTop}
\bysame, \emph{Analytic topology of groups, actions, strings and varieties},
  Geometry and dynamics of groups and spaces, Progr. Math., vol. 265,
  Birkh$\ddot{\text{a}}$user, Basel, 2008, pp.~3–--93.

\bibitem{Ruelle}
D.~Ruelle, \emph{Rotation numbers for diffeomorphisms and flows}, Ann. Inst. H.
  Poincar$\acute{\text{e}}$ Phys. Th$\acute{\text{e}}$or. \textbf{42} (1985),
  no.~1, 109–--115.

\bibitem{LoopRemarks}
E.~Shelukhin, \emph{Remarks on invariants of {H}amiltonian loops}, J. Topol.
  Anal. \textbf{2} (2010), no.~3, 277--325.

\bibitem{ShternAutomaticContinuity}
A.~I. Shtern, \emph{Automatic continuity of pseudocharacters on semisimple
  {L}ie groups}, Mat. Zametki \textbf{80} (2006), no.~3, 456--464, translation
  in Math. Notes 80 (2006), no. 3-4, 435–441.

\bibitem{Siegel}
C.~L. Siegel, \emph{Symplectic geometry}, Amer. J. Math. \textbf{65} (1943),
  no.~1, 1–--86.

\bibitem{GabiTobias2}
G.~Ben Simon and T.~Hartnick, \emph{Invariant orders on {H}ermitian {L}ie
  groups}, preprint arXiv:1011.3505v1, 2010.

\bibitem{GabiTobias3}
\bysame, \emph{Quasi total orders and translation numbers}, preprint
  www.math.ethz.ch/~bgabi/, 2010.

\bibitem{GabiTobias1}
\bysame, \emph{Reconstructing quasimorphisms from associated partial orders and
  a question of {P}olterovich}, preprint arXiv:0811.2608v4, 2010.

\bibitem{TosattiWeinkoveYauTamingSymplecticForms}
V.~Tosatti, B.~Weinkove, and S.-T. Yau, \emph{Taming symplectic forms and the
  {C}alabi-{Y}au equation}, Proc. Lond. Math. Soc. (3) \textbf{97} (2008),
  no.~2, 401--–424.

\bibitem{UhlenbeckYauGeneralHitchinKobayashi}
K.~Uhlenbeck and S.-T. Yau, \emph{On the existence of
  {H}ermitian-{Y}ang-{M}ills connections in stable vector bundles}, Comm. Pure
  Appl. Math. \textbf{39} (1986), no.~S, suppl., S257–--S293, Special Issue.
  Frontiers of the mathematical sciences: 1985 (New York, 1985).

\bibitem{UsherSpectralNumbersFloer}
M.~Usher, \emph{Spectral numbers in {F}loer theories}, Compos. Math.
  \textbf{144} (2008), no.~6, 1581--1592.

\bibitem{UsherDeformedQHandCalabiQm}
\bysame, \emph{Deformed {H}amiltonian {F}loer theory, capacity estimates, and
  {C}alabi quasimorphisms}, preprint arXiv:1006.5390v1, 2010.

\bibitem{UsherDualityFloerNovikovFiltration}
\bysame, \emph{Duality in filtered {F}loer-{N}ovikov complexes}, J. Topol.
  Anal. \textbf{2} (2010), no.~2, 233--258.

\bibitem{Weinstein}
A.~Weinstein, \emph{Cohomology of symplectomorphism groups and critical values
  of {H}amiltonians}, Math. Z. \textbf{201} (1989), no.~1, 75--82.

\end{thebibliography}
\end{document}